\input amstex
\documentstyle{amsppt}
%
%
\nologo
\nopagenumbers
\TagsOnRight
\def\const{\operatorname{const}}
\def\negskp{\hskip -2pt}
\def\compos{\,\raise 1pt\hbox{$\sssize\circ$} \,}
\pagewidth{360pt}
\pageheight{606pt}
\topmatter
\title Normal shift in general Lagrangian dynamics.
\endtitle
\author
Ruslan A. Sharipov
\endauthor
\abstract
It is well known that Lagrangian dynamical systems naturally
arise in describing wave front dynamics in the limit of short
waves (which is called pseudoclassical limit or limit of
geometrical optics). Wave fronts are the surfaces of constant
phase, their points move along lines which are called rays.
In non-homogeneous anisotropic media rays are not straight
lines. Their shape is determined by modified Lagrange equations.
An important observation is that for most usual cases propagating
wave fronts are perpendicular to rays in the sense of some
Riemannian metric. This happens when Lagrange function is quadratic
with respect to components of velocity vector. The goal of paper is
to study how this property transforms for the case of general
(non-quadratic) Lagrange function.
\endabstract
\address Rabochaya street 5, 450003, Ufa, Russia
\endaddress
\email \vtop to 20pt{\hsize=280pt\noindent
R\_\hskip 1pt Sharipov\@ic.bashedu.ru\newline
r-sharipov\@mail.ru\vss}
\endemail
\urladdr
http:/\negskp/www.geocities.com/r-sharipov
\endurladdr
\subjclass 53B40, 70H03, 70H05
\endsubjclass
\keywords
Wave fronts, normal shift, Lagrangian dynamics
\endkeywords
\endtopmatter
\loadbold
\document
\head
1. A simple preliminary example.
\endhead
    Description of most wave phenomena is based on wave equation.
This is second order partial differential equation of the following
form:
$$
\hskip -2em
\frac{1}{c^2}\,\frac{\partial^2\psi}{\partial {t^{\vphantom{i}}}^2}
-\sum^3_{i=1}\frac{\partial^2\psi}{\partial {x^i}^2}=0.
\tag1.1
$$
Here $t$ is time variable, while $x^1$, $x^2$, and $x^3$ are spatial
Cartesian coordinates. Parameter $c$ in first term is the velocity
of wave process described by the equation \thetag{1.1}. This is sound
velocity for sound waves in gases, liquids, or solid materials, and
this is light velocity for light waves in refracting media. For
homogeneous media $c$ is constant, but below we consider non-homogeneous
media, where $c=c(t,x^1,x^2,x^3)$.
\par
    Function $\psi=\exp(i\,(\omega\,t-k_1\,x^1-k_2\,x^2-k_3\,x^3))$ is
a solution of wave equation \thetag{1.1} for the case $c=\const$. It
describes a plane wave. Here $\omega$ is a frequency of wave, while
$k_1$, $k_2$, and $k_3$ are components of wave vector $\bold k$.
Frequency $\omega$ and wave vector $\bold k$ are related with each
other as follows:
$$
\hskip -2em
\omega=c\cdot|\bold k|.
\tag1.2
$$
The relationship \thetag{1.2} is called {\it dispersion law}. Short
wave limit corresponds to the case of high frequency, when $\omega\to
\infty$. Below we consider this case for non-homogeneous media with
$c\neq\const$. Therefore we cannot use simple exponential solution
$\psi=\exp(i\,(\omega\,t-k_1\,x^1-k_2\,x^2-k_3\,x^3))$ of wave
equation \thetag{1.1}. However, we can look for the exponential
solution with large parameter $\lambda\to\infty$:
$$
\hskip -2em
\psi=\sum^\infty_{\alpha=0}\frac{\psi_{\sssize(\ssize\alpha\sssize)}}
{(i\,\lambda)^\alpha}\cdot e^{i\lambda S}.
\tag1.3
$$
Substituting \thetag{1.3} into the equation \thetag{1.1}, we get the
following equation for $S$:
$$
\hskip -2em
\frac{1}{c^2}\left(\frac{\partial S}{\partial t}\right)^2
-\sum^3_{i=1}\left(\frac{\partial S}{\partial x^i}\right)^2=0.
\tag1.4
$$
This is well-known {\it eikonal equation} (see
Chapter~\uppercase\expandafter{\romannumeral 7} in \cite{1}). Suppose
that refracting properties of medium do not change in time. Then
$c=c(x^1,x^2,x^3)$. In this case we can consider a wave with constant
frequency $\omega=\lambda$. For such wave function $S$ in eikonal
equation \thetag{1.4} is taken to be linear function in time variable
$t$:
$$
\hskip -2em
S=t-\varphi(x^1,x^2,x^3).
\tag1.5
$$
Eikonal equation \thetag{1.4} then is written as the equation
for gradient of $\varphi$:
$$
\hskip -2em
\sum^3_{i=1}\left(\frac{\partial\varphi}{\partial x^i}\right)^2
-\frac{1}{c^2}=0.
\tag1.6
$$
Let's denote $-1/(2c^2)=U$ and then let's write the equation
\thetag{1.6} as
$$
\hskip -2em
\sum^3_{i=1}\frac{(\nabla_i\varphi)^2}{2}+U(x^1,x^2,x^3)=0.
\tag1.7
$$
Here $\nabla_1\varphi$, $\nabla_2\varphi$, and $\nabla_3\varphi$ are
components of gradient $\nabla\varphi$. We denote it by $\bold p$
and treat as a vector field in tree-dimensional space $\Bbb R^3$:
$$
\hskip -2em
\bold p=\nabla\varphi
=\Vmatrix
\partial\varphi/\partial x^1\\
\vspace{2ex}
\partial\varphi/\partial x^2\\
\vspace{2ex}
\partial\varphi/\partial x^3
\endVmatrix
\tag1.8
$$
If we substitute components of vector \thetag{1.8} into \thetag{1.7},
then we can write \thetag{1.7} as
$$
\hskip -2em
H(\nabla_1\varphi,\nabla_2\varphi,\nabla_3\varphi,x^1,x^2,x^3)=0,
\tag1.9
$$
where function $H=H(p_1,p_2,p_3,x^1,x^2,x^3)$ looks like Hamilton
function of a particle of unit mass $m=1$ in potential field
$U=U(x^1,x^2,x^3)$:
$$
\hskip -2em
H=\sum^3_{i=1}\frac{(p_i)^2}{2}+U.
\tag1.10
$$
\par
    Let $\varphi=\varphi(x^1,x^2,x^3)$ be a solution of the equation
\thetag{1.7} and let $\bold p=\nabla\varphi$ be corresponding momentum
vector field \thetag{1.8}. Let's consider integral curves of vector
field $\bold p$. They form two-parametric family of curves in
$\Bbb R^3$
$$
\hskip -2em
\cases
x^1=x^1(t,y^1,y^2),\\
x^2=x^2(t,y^1,y^2),\\
x^3=x^3(t,y^1,y^2)
\endcases
\tag1.11
$$
defined by solutions of the following system of ordinary differential
equations:
$$
\xalignat 3
&\hskip -2em
\dot x^1=p_1,&&\dot x^2=p_2,&&\dot x^3=p_3.
\tag1.12
\endxalignat
$$
Differential equations \thetag{1.12} can be written as Hamilton
equations:
$$
\hskip -2em
\dot x^i=\frac{\partial H}{\partial p_i},\qquad i=1,\,2,\,3.
\tag1.13
$$
Now let's calculate time derivative $\dot\bold p$ for the momentum
vector \thetag{1.8} due to the dynamics determined by differential
equations \thetag{1.12}:
$$
\dot p_i=\sum^3_{k=1}\frac{\partial p_i}{\partial x^k}\cdot\dot x^k=
\sum^3_{k=1}\frac{\partial^2\varphi}{\partial x^i\,\partial x^k}\cdot
\dot x^k=\sum^3_{k=1}\frac{\partial H}{\partial p_k}\cdot
\frac{\partial^2\varphi}{\partial x^i\,\partial x^k}.
$$
Remember that $\varphi$ is a solution of the equation \thetag{1.9}.
Differentiating \thetag{1.9}, we get
$$
\frac{\partial H(\nabla_1\varphi,\nabla_2\varphi,\nabla_3\varphi,
x^1,x^2,x^3)}{\partial x^i}=\sum^3_{k=1}\frac{\partial H}{\partial p_k}
\cdot\frac{\partial^2\varphi}{\partial x^i\,\partial x^k}
+\frac{\partial H}{\partial x^i}=0.
$$
Comparing the above two equalities, we derive differential equations
$$
\hskip -2em
\dot p_i=-\frac{\partial H}{\partial x^i},\qquad i=1,\,2,\,3.
\tag1.14
$$
Both \thetag{1.13} and \thetag{1.14} form complete system of Hamilton
equations 
$$
\xalignat 2
&\hskip -2em
\dot x^i=\frac{\partial H}{\partial p_i},
&&\dot p_i=-\frac{\partial H}{\partial x^i}
\tag1.15
\endxalignat
$$
with Hamilton function \thetag{1.10}.\par
    Note that Hamilton equations \thetag{1.15} is a system of 6 first
order ODE's. Its solutions define five-parametric family of curves
in $\Bbb R^3$:
$$
\hskip -2em
\cases
x^1=x^1(t,y^1,y^2,y^3,y^4,y^5),\\
x^2=x^2(t,y^1,y^2,y^3,y^4,y^5),\\
x^3=x^3(t,y^1,y^2,y^3,y^4,y^5).
\endcases
\tag1.16
$$
Curves \thetag{1.11} form two-parametric subfamily in five-parametric
family of curves \thetag{1.16}. They are distinguished by the following
two properties:
\rosteritemwd=4pt
\roster
\item"1)" curves \thetag{1.11} correspond to zero level of energy $H=0$;
\item"2)" curves \thetag{1.11} are perpendicular to level surfaces of
the function $\varphi(x^1,x^2,x^3)$.
\endroster
First property follows from \thetag{1.9}. Second is obvious, since curves
\thetag{1.11} are directed along gradient vector \thetag{1.8}. One can
calculate complete derivative of the function $\varphi(x^1,x^2,x^3)$ with
respect to parameter $t$ along these curves:
$$
\hskip -2em
\frac{d\varphi}{dt}=\Omega=\sum^3_{i=1}p_i\,\frac{\partial H}{\partial p_i}.
\tag1.17
$$
Curves \thetag{1.16} defined by Hamilton equations \thetag{1.15} and
restricted by the above conditions 1) and 2) are called {\it characteristic
lines} for nonlinear first order partial differential equation \thetag{1.9}.
They are used in order to construct solutions of this equation as described
just below (see also \cite{2} and \cite{3}).\par
    Let's take some smooth surface $\sigma$ in $\Bbb R^3$. We assume that
$\sigma$ is level surface with $\varphi=0$ for the solution
$\varphi(x^1,x^2,x^3)$ of the equation \thetag{1.9} that we are going to
construct. Denote by $y^1$ and $y^2$ inner curvilinear coordinates of points
on $\sigma$. Then we can write the equations determining points of
$\sigma$ in parametric form:
$$
\hskip -2em
\cases
x^1=x^1(y^1,y^2),\\
x^2=x^2(y^1,y^2),\\
x^3=x^3(y^1,y^2).
\endcases
\tag1.18
$$
At each point of $\sigma$ we have unit normal vector $\bold n$. Let's denote
it by $\bold n=\bold n(y^1,y^2)$. Assuming $\sigma$ to be orientable, we can
take $\bold n(y^1,y^2)$ to be smooth function of $y^1$ and $y^2$. Under these
assumptions we define vector function
$$
\hskip -2em
\bold p=\nu\cdot\bold n
\tag1.19
$$
on $\sigma$, getting scalar factor $\nu=\nu(y^1,y^2)$ from the following
equality:
$$
\hskip -2em
H(\nu\,n_1,\,\nu\,\,n_2,\,\nu\,n_3,x^1,x^2,x^3)=0.
\tag1.20
$$
Then we use vector function \thetag{1.19} in order to set up Cauchy problem
$$
\xalignat 2
&\hskip -2em
x^i\,\hbox{\vrule height 8pt depth 8pt width 0.5pt}_{\,t=0}
=x^i(y^1,y^2),
&&p_i\,\hbox{\vrule height 8pt depth 8pt width 0.5pt}_{\,t=0}=
p_i(y^1,y^2)
\tag1.21
\endxalignat
$$
for Hamilton equations \thetag{1.15}. Solving this Cauchy problem
\thetag{1.21}, we
obtain two-parametric family of characteristic lines given by functions
\thetag{1.11} that extend initial functions \thetag{1.18}. They possess
property 1), since we determine $\nu$ by \thetag{1.20}. They also possess
property 2), since we determine $\bold p$ by \thetag{1.19} (at least for
initial surface $\sigma$). These characteristic lines fill some neighborhood
of initial surface $\sigma$. Therefore we can treat $t$, $y^1$, $y^2$ as
curvilinear coordinates in $\Bbb R^3$ and consider \thetag{1.11} as
transition functions to these curvilinear coordinates. Then integral
$$
\hskip -2em
\varphi=\int\limits^{\,t}_{\!0}\Omega\,dt,
\tag1.22
$$
where $\Omega$ is given by right hand side of \thetag{1.17}, yields a
solution of partial differential equation \thetag{1.9} expressed in
curvilinear coordinates $t$, $y^1$, and $y^2$. This solution \thetag{1.22}
satisfies zero boundary-value condition on $\sigma$:
$$
\hskip -2em
\varphi\,\hbox{\vrule height 8pt depth 8pt width 0.5pt}_{\,\sigma}=0.
\tag1.23
$$
In other words, \thetag{1.23} means that $\sigma$ is zero level surface
for the function $\varphi$.\par
    Note that in curvilinear coordinates $t$, $y^1$, $y^2$ initial surface
$\sigma$ is given by the equation $t=0$. However, other level surfaces of
the function $\varphi$ are not given by the equations $t=\const$. In order
to change this situation we should choose another set of curvilinear
coordinates $s$, $y^1$, $y^2$, where $s=\varphi(t,y^1,y^2)$. This means,
that we change parametrization of characteristic lines \thetag{1.11}
without changing them as geometric sets of points. In new parameter $s$
characteristic lines of the equation \thetag{1.9} are given by modified
Hamilton equations
$$
\xalignat 2
&\hskip -2em
\dot x^i=\frac{1}{\Omega}\,\frac{\partial H}{\partial p_i},
&&\dot p_i=-\frac{1}{\Omega}\,\frac{\partial H}{\partial x^i},
\tag1.24
\endxalignat
$$
where denominator $\Omega$ is determined by right hand side of
\thetag{1.17}. In our particular case, when function $H$ is given
by formula \thetag{1.10}, we have $\Omega=(p_1)^2+(p_2)^2
+(p_3)^2$. Hence $\Omega\neq 0$ for $\bold p\neq 0$.\par
    Let's fix new curvilinear coordinates $s$, $y^1$, $y^2$. Here
$\varphi(s,y^1,y^2)=s$ by definition. Now let's return to initial
wave equation \thetag{1.1} and to formula \thetag{1.5} for the
function $S$ in asymptotical power expansion \thetag{1.3}. It's
important to note that $t$ in \thetag{1.5} do not coincide with $t$
in \thetag{1.11} and in Hamilton equations \thetag{1.15}, where $t$
was used as a parameter on characteristic lines of the equation
\thetag{1.6}. Therefore now in the expression for $S$ we have both
$s$ and $t$ (and $t$ is time variable again):
$$
S=S(t,s,y^1,y^2)=t-s.
$$
For exponential factor $e^{i\lambda S}$ in \thetag{1.3}, taking into
account that $\lambda=\omega$, we get:
$$
\hskip -2em
e^{i\lambda S}=e^{i\omega(t-s)}.
\tag1.25
$$
Right hand side of \thetag{1.25} corresponds to plane wave propagating
in the direction of $s$-axis. In original Cartesian coordinates $x^1,\,
x^2,\,x^3$ this looks like non-plain wave propagating along characteristic
lines of the equation \thetag{1.9}. Level surfaces of the function
$\varphi$ are the surfaces of constant phase in such wave. They are called
{\it wave fronts}. The equation $t-s=\const$, when transformed to Cartesian
coordinates $x^1,\,x^2,\,x^3$, describes moving surface, that gradually
passes positions of level surfaces of the function $\varphi$. This process
is called wave front dynamics. It's very important that this process can
be understood as a {\bf motion of separate points of wave front}, each
obeying modified Hamilton equations \thetag{1.24}. For this reason these
equations are called {\bf the equations of wave front dynamics}.\par
    Another important point concerning wave front dynamics, that we noted
above, is that level surfaces of the function $\varphi$ are perpendicular
to characteristic lines \thetag{1.11}. Therefore wave front dynamics is
a normal displacement (or {\bf normal shift}) of initial surface $\sigma$
along trajectories of modified Hamiltonian dynamical system.
\head
2. More complicated example.
\endhead
    Let $M$ be some Riemannian manifold. Denote by $\nabla$ standard
covariant differentiation determined by metric connection $\Gamma$ in
$M$. The following differential operator is called Laplace-Beltrami
operator in the manifold $M$:
$$
\hskip -2em
\triangle=\sum^n_{i=1}\sum^n_{j=1}g^{ij}\,\nabla_i\,\nabla_j.
\tag2.1
$$
Here $g^{ij}$ are components of metric tensor in local coordinates
$x^1,\,\ldots,\,x^n$, while $\nabla_i$ and $\nabla_j$ are symbols of
covariant derivatives in these local coordinates. Differential operator
$H$ is called {\it fiberwise spherically symmetric} if it is represented
as a polynomial of Laplace-Beltrami operator \thetag{2.1}:
$$
\hskip -2em
H(p,D)=\sum^m_{k=0}a_k(p)\,\triangle^k.
\tag2.2
$$
Here $p$ is a point of $M$ and $D$ is a formal symbol for differentiation.
Coefficients $a_0,\,\ldots,\,a_m$ in \thetag{2.2} are arbitrary smooth
functions of $p\in M$. Note that $H$ is scalar operator. Operator
\thetag{2.2} can be applied either to scalar field or tensorial field
in $M$, yielding the field of the same type as that it was applied to.
\par
    Now let's add differentiation in time variable $\partial_t=\partial/
\partial t$ and let's introduce large parameter $\lambda$ to \thetag{2.2}.
As a result we get differential operator
$$
H(p,\lambda^{-1}D)=\sum^m_{s=0}\sum^m_{k=0}\frac{a_{sk}(p)}
{(i\,\lambda)^{s+2k}}\,{\partial_t}^s\,\triangle^k.
$$
The following differential equation in $M$ is an analog of wave
equation \thetag{1.1}:
$$
\hskip -2em
H(p,\lambda^{-1}D)\psi=0.
\tag2.3
$$
Short wave asymptotics $\lambda\to\infty$ for this equations is described
by the same asymptotical expansion \thetag{1.3} as in case of standard
wave equation. Coefficients $a_{sk}(p)$ in \thetag{2.3} do not depend on
$t$. Therefore we can choose $S$ to be linear function of $t$:
$$
\hskip -2em
S=t-\varphi(p),
\tag2.4
$$
just like it was in \thetag{1.5}. Substituting \thetag{2.4} into
\thetag{1.3} and substituting \thetag{1.3} into \thetag{2.3}, we
derive differential equation for phase function $\varphi(p)$ in
\thetag{2.4}:
$$
\hskip -2em
\sum^m_{k=0}\left(\,\shave{\sum^m_{s=0}}a_{sk}(p)\right)
\cdot|\nabla\varphi|^{2k}=\sum^m_{k=0}b_k(p)\cdot|\nabla
\varphi|^{2k}=0.
\tag2.5
$$
Here $|\nabla\varphi|$ is modulus of covector field $\nabla\varphi$
measured in Riemannian metric $\bold g$:
$$
|\nabla\varphi|^2=\sum^n_{i=1}\sum^n_{j=1}g^{ij}\,\nabla_i\varphi\,
\nabla_j\varphi.
$$
Let's denote $\nabla\varphi$ by $\bold p$ as it was done above
in section~1 (see formula \thetag{1.8}):
$$
\bold p=\nabla\varphi
=\Vmatrix
\partial\varphi/\partial x^1\\
\vspace{2ex}
\vdots
\vspace{2ex}
\partial\varphi/\partial x^n
\endVmatrix.
$$
Now we can write \thetag{2.5} as polynomial equation with respect
to components of $\bold p$:
$$
\hskip -2em
\sum^m_{k=0}b_k(p)\cdot|\bold p|^{2k}=0.
\tag2.6
$$
Let's denote left hand side of \thetag{2.6} by $H=H(p,\bold p)$.
The equation \thetag{2.5}, written as
$$
\hskip -2em
H(p,\nabla\varphi)=0,
\tag2.7
$$
is exact analog of the equation \thetag{1.9} from section~1. Further
steps in solving this equation are quite similar to those in section~1
(they are described in details in paper \cite{4}). Below we shall
not discuss them. However, we shall point out most important features
of wave front dynamics in the limit of short waves $\lambda=\omega\to
\infty$ for generalized wave equation \thetag{2.3}. They are the
following ones:
\roster
\rosteritemwd=1pt
\item"---" wave fronts are level hypersurfaces $\sigma_t=\{p\in M:\
\varphi(p)=t\}$ for the function $\varphi(p)$, where $\varphi(p)$
is a solution of differential equation \thetag{2.7};
\item"---" time evolution of wave fronts $\sigma_t$ in $M$ can be
described in terms of motion of their points obeying modified Hamilton
equations
$$
\xalignat 2
&\hskip -2em
\dot x^i=\frac{1}{\Omega}\,\frac{\partial H}{\partial p_i},
&&\dot p_i=-\frac{1}{\Omega}\,\frac{\partial H}{\partial x^i},
\tag2.8
\endxalignat
$$
where Hamilton function $H$ is determined by left hand side of \thetag{2.6}
$$
\hskip -2em
H(p,\bold p)=\sum^m_{k=0}b_k(p)\cdot|\bold p|^{2k},
\tag2.9
$$
and denominator $\Omega$ in \thetag{2.8} is determined by formula
$$
\Omega=\sum^n_{i=1}p_i\,\frac{\partial H}{\partial p_i};
$$
\item"---" wave front dynamics for wave equation \thetag{2.3} in
short wave limit $\lambda\to\infty$ is a normal shift of of initial
wave front hypersurface $\sigma$ along trajectories of modified
Hamiltonian dynamical system \thetag{2.8}, this means that
orthogonality of wave fronts $\sigma_t$ and trajectories of shift
is preserved in time;
\item"---" normal shift of hypersurface $\sigma$ is initiated by Cauchy
problem data
$$
\xalignat 2
&\hskip -2em
x^i\,\hbox{\vrule height 8pt depth 8pt width 0.5pt}_{\,t=0}
=x^i(p),
&&p_i\,\hbox{\vrule height 8pt depth 8pt width 0.5pt}_{\,t=0}=
\nu(p)\cdot n_i(p)
\tag2.10
\endxalignat
$$
for the equations \thetag{2.8}, where $n_i(p)$ are covariant components
of normal vector $\bold n(p)$ at the point $p\in\sigma$, and $\nu(p)$
is a scalar factor determined by the equation
$$
\hskip -2em
H(p,\nu\cdot\bold n(p))=0.
\tag2.11
$$
\endroster
   For us the most important feature of wave front dynamics, among those
listed above, is the phenomenon of normal shift. It was revealed in simplest
case considered in section~1. It is also present in more complicated case
related to some Riemannian metric. Our aim below is to reveal this
phenomenon for the case of general Hamilton function $H$, which is not
restricted by formula \thetag{2.9}. In order to do this we need to
introduce geometrical technique, which is not new, but nevertheless, is
not commonly known. It seems to me, that this technique first appeared in
Finslerian geometry (see \cite{5} and \cite{6}). We used this technique
in \cite{7--22}, where theory of Newtonian dynamical systems admitting
normal shift was developed (see also theses [23], [24], and recent papers
\cite{25--29}).
\head
3. Extended tensor fields.
\endhead
    Let's consider Hamilton function \thetag{2.9}. It depends on two
arguments $p$ and $\bold p$, where $p$ is a point of manifold $M$, while
$\bold p$ is cotangent vector at the point $p$, i\.\,e\. $\bold p$ is
an element of cotangent space $T^*_p(M)$. Both $p$ and $\bold p$, taken
together, form a pair $q=(p,\bold p)$ which is a point of cotangent bundle
$T^*\!M$. This means that $H$ is a scalar field in cotangent bundle
$T^*\!M$. But we shall treat it as {\it extended scalar field} in $M$
as defined below. Let's consider the following tensor product:
$$
T^r_s(p,M)=\overbrace{T_p(M)\otimes\ldots\otimes T_p(M)}^{\text{$r$
times}}\otimes\underbrace{T^*_p(M)\otimes\ldots\otimes T^*_p(M)}_{\text{$s$
times}}
$$
Tensor product $T^r_s(p,M)$ is known as a space of $(r,s)$-tensors at the
point $p\in M$. Pair of integer numbers $(r,s)$ determines type of tensors.
Elements of $T^r_s(p,M)$ are called {\it $r$-times contravariant and
$s$-times covariant tensors} or simply $(r,s)$-tensors.
\definition{Definition 3.1} {\it Extended} tensor field $\bold X$
of type $(r,s)$ in $M$ is a tensor-valued function that maps each
point $q=(p,\bold p)$ of some domain $G\subseteq T^*\!M$ to a tensor
of the space $T^r_s(p,M)$. If $G=T^*\!M$, then $\bold X$ is called
{\it global} extended tensor field.
\enddefinition
   Note a trick: arguments of extended tensor fields belong to cotangent
bundle $T^*\!M$, while their values are tensors related to base manifold
$M$. If we replace $T^*\!M$ by tangent bundle $TM$, we can state another
definition of extended tensor field.
\definition{Definition 3.2} {\it Extended} tensor field $\bold X$
of type $(r,s)$ in $M$ is a tensor-valued function that maps each
point $q=(p,\bold p)$ of some domain $G\subseteq TM$ to a tensor
of the space $T^r_s(p,M)$. If $G=TM$, then $\bold X$ is called
{\it global} extended tensor field.
\enddefinition
    In the case of arbitrary smooth manifold $M$ definitions~3.1 and
3.2 lead to different theories. But for Riemannian manifold $M$ tangent
bundle $TM$ and cotangent bundle $T^*\!M$ are bound with each other
by duality maps:
$$
\xalignat 2
&\hskip -2em
\bold g\!:TM\to T^*\!M,
&&\bold g^{-1}\!:T^*\!M\to TM.
\tag3.1
\endxalignat
$$
In local coordinates duality maps \thetag{3.1} are represented as
index lowering and index raising procedures in arguments of
extended tensor field $\bold X$:
$$
\xalignat 2
&p_i=\sum^n_{j=1}g_{ij}\,p^j,
&&p^i=\sum^n_{j=1}g^{ij}\,p_j.
\endxalignat
$$
Due to duality maps \thetag{3.1} two objects introduced by
definitions~3.1 and 3.2 are the same in essential. We call them
{\it covariant} and {\it contravariant} representations of extended
tensor field $\bold X$.\par
    In local coordinates $x^1,\,\ldots,\,x^n$ extended tensor field
$\bold X$ is represented by its components $X^{i_1\ldots\,i_r}_{j_1
\ldots\,j_s}(x^1,\ldots,x^n,p_1,\ldots,p_n)$ or $X^{i_1\ldots\,
i_r}_{j_1\ldots\,j_s}(x^1,\ldots,x^n,p^1,\ldots,p^n)$, depending on
which representation (covariant or contravariant) is used. Extended
tensor field $\bold X$ is called {\it smooth} if its components are
smooth functions.\par
     Smooth extended tensor fields form a ring, we denote it by $\goth
F=\goth F(T^*\!M)$ in the case of covariant representation, and by
$\goth F=\goth F(TM)$ in the case of contravariant representation.
The whole set of smooth extended tensor fields in $M$ is equipped
with operations of 1) summation, 2) multiplications by scalars, 3)
tensor product, 4) contraction. It forms bi-graded algebra over the
ring $\goth F$. We denote this algebra by $\bold T(M)$ and call it
{\it an algebra of extended tensor fields} in $M$:
$$
\hskip -2em
\bold T(M)=\bigoplus^\infty_{r=0}\bigoplus^\infty_{s=0}T^r_s(M).
\tag3.2
$$
\definition{Definition 3.3} A map $D\!:\bold T(M)\to\bold T(M)$ is
called a {\it differentiation} of extended algebra of tensor fields,
if the following conditions are fulfilled:
\roster
\rosteritemwd=10pt
\item concordance with grading: $D(T^r_s(M))\subset T^r_s(M)$;
\item $\Bbb R$-linearity: $D(\bold X+\bold Y)=D(\bold X)+D(\bold Y)$
      and $D(\lambda\bold X)=\lambda D(\bold X)$ for $\lambda\in\Bbb R$;
\item commutation with contractions: $D(C(\bold X))=
      C(D(\bold X))$;
\item Leibniz rule: $D(\bold X\otimes\bold Y)=D(\bold X)
      \otimes\bold Y+\bold X\otimes D(\bold Y)$.
\endroster
\enddefinition
\noindent Theory of differentiations in extended algebra of tensor
fields \thetag{3.2} is considered
in Chapters~\uppercase\expandafter{\romannumeral 2},
\uppercase\expandafter{\romannumeral 3}, and \uppercase
\expandafter{\romannumeral 4} of thesis \cite{23}. In this section
below we shall mention some facts from this theory needed for
further use.\par
    Suppose that $\bold T(M)$ is extended algebra of tensor fields
in $M$ taken in contravariant representation. Then the set of its
differentiations $\goth D(M)$ possesses the structure of module over
the ring $\goth F(TM)$. The set of extended vector fields (i\.\,e\.
summand $T^1_0(M)$ in direct sum \thetag{3.2}) also possesses the structure
of $\goth F(TM)$-module. Therefore the following definition is consistent.
\definition{Definition 2.1} {\it Covariant differentiation} $\nabla$ in the
algebra of extended tensor fields $\bold T(M)$ is a homomorphism of $\goth
F(TM)$-modules $\nabla\!:T^1_0(M)\to\goth D(M)$. Image of vector field
$\bold Y$ under such homomorphism denoted by $\nabla_{\bold Y}$ is called
{\it covariant differentiation along vector field $Y$}.
\enddefinition
    For each covariant differentiation the expression $\nabla_{\bold Y}
\bold X$ is $\goth F(TM)$-linear with respect to $\bold Y$. Therefore
$\nabla$ can be treated as a map $\nabla\!:T^r_s(M)\to T^r_{s+1}(M)$.
Each smooth manifold $M$ possesses exactly one canonical covariant
differentiation $\tilde\nabla$ which is called {\it vertical gradient}.
In local coordinates it is expressed by formula
$$
\hskip -2em
\tilde\nabla_qX^{i_1\ldots\,i_r}_{j_1\ldots\,j_s}=\frac{\partial
X^{i_1\ldots\,i_r}_{j_1\ldots\,j_s}}{\partial p^q}.
\tag3.3
$$
In order to define other covariant differentiations one need some
additional geometric structures in $M$. Thus, if $M$ possesses affine
connection $\Gamma$, one can define {\it horizontal gradient} $\nabla$.
In local coordinates it is expressed by formula
$$
\hskip -2em
\aligned
&\nabla_qX^{i_1\ldots\,i_r}_{j_1\ldots\,j_s}=\frac{\partial
X^{i_1\ldots\,i_r}_{j_1\ldots\,j_s}}{\partial x^q}
-\sum^n_{a=1}\sum^n_{b=1}p^a\,\Gamma^b_{qa}\,\frac{\partial
X^{i_1\ldots\,i_r}_{j_1\ldots\,j_s}}{\partial p^b}\,+\\
&+\sum^r_{k=1}\sum^n_{a_k=1}\!\Gamma^{i_k}_{q\,a_k}\,X^{i_1\ldots\,
a_k\ldots\,i_r}_{j_1\ldots\,\ldots\,\ldots\,j_s}
-\sum^s_{k=1}\sum^n_{b_k=1}\!\Gamma^{b_k}_{q\,j_k}
X^{i_1\ldots\,\ldots\,\ldots\,i_r}_{j_1\ldots\,b_k\ldots\,j_s}.
\endaligned
\tag3.4
$$
If we take covariant representation of the algebra of extended
tensor fields $\bold T(M)$, then formulas for vertical and horizontal
gradients are transformed as follows:
$$
\align
&\hskip -2em
\tilde\nabla^qX^{i_1\ldots\,i_r}_{j_1\ldots\,j_s}=
\frac{\partial X^{i_1\ldots\,i_r}_{j_1\ldots\,j_s}}
{\partial p_q}.
\tag3.5\\
\vspace{2ex}
&\hskip -2em
\aligned
&\nabla_{\!q}X^{i_1\ldots\,i_r}_{j_1\ldots\,j_s}=\frac{\partial
X^{i_1\ldots\,i_r}_{j_1\ldots\,j_s}}{\partial x^q}
+\sum^n_{a=1}\sum^n_{b=1}p_a\,\Gamma^a_{qb}\,\frac{\partial
X^{i_1\ldots\,i_r}_{j_1\ldots\,j_s}}{\partial p_b}\,+\\
&+\sum^r_{k=1}\sum^n_{a_k=1}\!\Gamma^{i_k}_{q\,a_k}\,X^{i_1\ldots\,
a_k\ldots\,i_r}_{j_1\ldots\,\ldots\,\ldots\,j_s}
-\sum^s_{k=1}\sum^n_{b_k=1}\!\Gamma^{b_k}_{q\,j_k}\,
X^{i_1\ldots\,\ldots\,\ldots\,i_r}_{j_1\ldots\,b_k\ldots\,j_s}.
\endaligned
\tag3.6
\endalign
$$
In the case of arbitrary smooth manifold $M$ gradients defined
by formulas \thetag{3.3} and \thetag{3.4} are not related to
those defined by formulas \thetag{3.5} and \thetag{3.6}. However,
if $M$ is Riemannian manifold, then $\tilde\nabla$ and $\nabla$
defined by these two ways appear to be the same
differentiations\footnotemark\ in different representations of
algebra $\bold T(M)$.\footnotetext{One should only lower index $q$
in \thetag{3.5}.}\adjustfootnotemark{-1} This fact is expressed by
the following commutation relationships:
$$
\xalignat 2
&\nabla(\bold X\compos\bold g)=(\nabla\bold X)\compos\bold g,
&&\tilde\nabla(\bold X\compos\bold g)=(\tilde\nabla\bold X)
\compos\bold g.
\endxalignat
$$
Here $\bold g$ is duality map \thetag{3.1} defined by metric
tensor of Riemannian manifold.
\head
4. Legendre transformation.
\endhead
   Legendre transformation is usually used to relate Lagrangian and
Hamiltonian dynamical systems. Suppose that $M$ is smooth manifold
and let $L(p,\bold v)$ be smooth extended scalar field in $M$ taken
in contravariant representation. Then dynamical system in tangent
bundle $TM$ described by differential equations
$$
\xalignat 2
&\hskip -2em
\dot x^i=v^i,
&&\frac{d}{dt}\left(\frac{\partial L}{\partial v^i}\right)=
\frac{\partial L}{\partial x^i}
\tag4.1
\endxalignat
$$
is called Lagrangian dynamical system. Let's apply covariant derivative
\thetag{3.5} to $L$. As a result we get covector field $\bold p=\bold p(p,
\bold v)$ with components
$$
\hskip -2em
p_i=\tilde\nabla_iL=\frac{\partial L}{\partial v^i}.
\tag4.2
$$
If pair $(p,\bold v)$ is a point of tangent bundle $TM$, then pair
$(p,\bold p)$ is a point of cotangent bundle $T^*\!M$. This means that
derivatives \thetag{4.2} determine a map
$$
\hskip -2em
\lambda\!:TM\to T^*\!M.
\tag4.3
$$
This map is known as Legendre transformation (see \cite{30}). Below
we assume Legendre transformation \thetag{4.3} to be invertible.
Moreover we assume inverse map
$$
\hskip -2em
\lambda^{-1}\!:T^*\!M\to TM
\tag4.4
$$
to be smooth. Under these assumptions we can treat direct and inverse
Legendre transformations \thetag{4.3} and \thetag{4.4} as nonlinear
analogs of duality maps \thetag{3.1}. Lagrange function $L(p,\bold v)$
and Lagrange equations \thetag{4.1} are associated with tangent bundle
$TM$. Vector $\bold v$ in arguments of Lagrange function is called
{\it velocity vector}, while covector $\bold p$ with components
\thetag{4.2} is called {\it momentum covector}. This gives rise to the
following terminology. If $\bold X$ is an extended tensor field
in contravariant representation and if $\bold Y=\bold X\compos
\lambda^{-1}$, then we say that $\bold Y$ is $\bold p$-representation
or {\it momentum representation} for $\bold X$, while $\bold X$ is
called $\bold v$-representation or {\it velocity representation} for
$\bold Y$. For the case of general smooth manifold $M$ (without
Riemannian metric) direct and inverse Legendre transformations
\thetag{4.3} and \thetag{4.4} bind the following two representations
of extended tensor fields:
$$
\vbox{\offinterlineskip
\def\vr{\vrule height 14pt depth 8pt}
\settabs\+\hskip 0.5cm&\hskip 4.5cm&\hskip 0.5cm&\hskip 5.2cm&\cr
\hrule
\+\vr\strut &covariant $\bold p$-representation
 &\vr\strut &contravariant $\bold v$-representation &\vr\strut\cr
\hrule}
$$
If $M$ is Riemannian manifold, we have four representations per each
extended field:
$$
\vbox{\offinterlineskip
\def\vr{\vrule height 14pt depth 8pt}
\settabs\+\hskip 0.5cm&\hskip 5.2cm&\hskip 0.5cm&\hskip 5.2cm&\cr
\hrule
\+\vr\strut &covariant $\bold p$-representation
 &\vr\strut &contravariant $\bold v$-representation &\vr\strut\cr
\hrule
\+\vr\strut &contravariant $\bold p$-representation
 &\vr\strut &covariant $\bold v$-representation &\vr\strut\cr
\hrule}
$$\par
    Now let's consider the following two extended scalar fields
$h$ and $H$:
$$
\xalignat 2
&\hskip -2em
h=\sum^n_{i=1}v^i\,\tilde\nabla_iL-L,
&&H=h\compos\lambda^{-1}.
\tag4.5
\endxalignat
$$
Scalar field $H$ is known as Hamilton function, while $h$ is its
$\bold v$-representation. Lagrange function $L$ and its
$\bold p$-representation $l=L\compos\lambda^{-1}$ can be expressed
through Hamilton function $H$ by formulas similar to \thetag{4.5}:
$$
\xalignat 2
&\hskip -2em
l=\sum^n_{i=1}p_i\,\tilde\nabla^iH-H,
&&L=l\compos\lambda.
\tag4.6
\endxalignat
$$
Applying Legendre transformation to Lagrangian dynamical system
\thetag{4.1}, we can transform them to Hamiltonian dynamical system
in cotangent bundle $T^*\!M$:
$$
\xalignat 2
&\dot x^i=\frac{\partial H}{\partial p_i},
&&\dot p_i=-\frac{\partial H}{\partial x^i}.
\endxalignat
$$
This fact is well known (see \cite{30}), as well as above formulas
\thetag{4.5} and \thetag{4.6}.
\head
5. Modified Hamiltonian and Lagrangian dynamical systems
in Riemannian manifolds.
\endhead
    Now suppose that manifold $M$ is equipped with Riemannian metric
$\bold g$. All results we discussed in section~2 were obtained under
this assumption. Remember that wave front dynamics is described by
modified Hamilton equations \thetag{2.8}. Using \thetag{3.5} and
\thetag{3.6}, we can replace partial derivatives in them by covariant
derivatives:
$$
\xalignat 2
&\hskip -2em
\dot x^i=\frac{\tilde\nabla^i H}{\Omega},
&&\nabla_tp_i=-\frac{\nabla_i H}{\Omega},
\tag5.1
\endxalignat
$$
Time derivatives $\dot x^1,\,\ldots,\,\dot x^n$ in \thetag{5.1} are
components of tangent vector to trajectory $p=p(t)$, while $\nabla_t$
is standard covariant derivative with respect to parameter $t$ along
this curve. Similar to original Hamilton equations \thetag{4.5},
modified Hamilton equations are associated with cotangent bundle
$T^*\!M$. Using inverse Legendre map \thetag{4.4}, one can transform them
to $\bold v$-representation. This was done in paper \cite{4}. As a
result modified Lagrange equations were obtained:
$$
\xalignat 2
&\hskip -2em
\dot x^i=\frac{v^i}{\Omega},
&&\nabla_t(\tilde\nabla_iL)=\frac{\nabla_iL}{\Omega}.
\tag5.2
\endxalignat
$$
Denominator $\Omega$ in original $\bold p$-representation is given
by formula
$$
\hskip -2em
\Omega=\sum^n_{i=1}p_i\,\tilde\nabla^iH,
\tag5.3
$$
Upon passing to $\bold v$-representation in \thetag{5.2} formula
\thetag{5.3} transforms to
$$
\hskip -2em
\Omega=\sum^n_{i=1}v^i\,\tilde\nabla_iL.
\tag5.4
$$
Formula \thetag{5.4} means that we can express $\Omega$ in
$\bold v$-representation explicitly through Lagrange function
$L$.
\head
6. Newtonian dynamical systems admitting normal shift
in Riemannian manifolds.
\endhead
    One can see that vector $\bold v$ for modified Lagrangian dynamics
\thetag{5.2} do not coincide with actual velocity vector. If we
denote actual velocity vector by $\bold u$, then we derive
$$
\hskip -2em
u^i=\frac{v^i}{\Omega}.
\tag6.1
$$
Formula \thetag{6.1} defines nonlinear map similar to Legendre map
$\lambda$:
$$
\hskip -2em
\mu\!: TM\to TM
\tag6.2
$$
If this map \thetag{6.2} is invertible and if inverse map $\mu^{-1}$ is
smooth, then one can transform modified Lagrange equations \thetag{5.2}
to the following form:
$$
\xalignat 2
&\hskip -2em
\dot x^i=u^i,
&&\nabla_tu^i=F^i(x^1,\ldots,x^n,u^1,\ldots,u^n).
\tag6.3
\endxalignat
$$
Differential equations \thetag{6.3} determine Newtonian dynamical system
in $M$. Extended vector field $\bold F$ with components $F^1,\,\ldots,\,
F^n$ in \thetag{6.3} is called {\it force field} of this Newtonian
dynamical system.\par
    Note that nonlinear map \thetag{6.2} is more complicated than Legendre
map. Almost each modified Lagrangian dynamical system can be transformed
to Newtonian form (at least locally). However, converse is not true.
Moreover, even if it is known that Newtonian dynamical system \thetag{6.3}
is derived from modified Lagrangian dynamical system \thetag{5.2}, there is
no explicit formula for $\bold F$. In paper \cite{4} one can find explicit
formula for $\bold F$, but in very special case, when Lagrange function
$L$ is {\it fiberwise spherically symmetric} with respect to Riemannian
metric $\bold g$:
$$
\hskip -2em
F_k=-|\bold u|\cdot\sum^n_{i=1}\frac{\nabla_iW}{W'}\cdot
\left(2\,N^i\,N_k-\delta^i_k\right).
\tag6.4
$$
Here $W=W(p,|\bold u|)=h\compos\mu^{-1}$ is $\bold u$-representation for
Hamilton function $H$ and $W'$ is partial derivative of $W$ with respect
to its $(n+1)$-th argument $u=|\bold u|$, $N^i$ and $N_k$ are components
of unit vector $\bold N=\bold u/|\bold u|$.\par
    Remember that in the example considered in section~2 Hamilton function
$H$ is given by formula \thetag{2.9}. It is fiberwise spherically symmetric,
i\.\,e\. $H=H(p,|\bold p|)$. Applying inverse Legendre map we get fiberwise
spherically symmetric function $h=h(p,|\bold v|)=H\compos\lambda$. Its
$\bold u$ representation then is fiberwise spherically symmetric function
$W$ in formula \thetag{6.4}. Thus we have the following theorem.
\proclaim{Theorem 6.1} For generalized wave equation \thetag{2.3} in
Riemannian manifold $M$ wave front dynamics in the limit of short
waves is described by Newtonian dynamical system \thetag{6.3} with force
field \thetag{6.4}.
\endproclaim
    Theorem~6.1 is the main result of paper \cite{4}. It establishes
a link between wave propagation phenomena and the {\it theory of dynamical
systems admitting normal shift} (see papers \cite{7--22}). Below we give
brief introduction to this theory.\par
    Let $M$ be Riemannian manifold and let $\sigma$ be some smooth
hypersurface in $M$. Suppose that $p$ is a point of $\sigma$ and
$\bold n(p)$ is a unit normal vector to $\sigma$ at the point $p$.
Under these assumptions we can consider initial data
$$
\xalignat 2
&\hskip -2em
x^i\,\hbox{\vrule height 8pt depth 8pt width 0.5pt}_{\,t=0}
=x^i(p),
&&u^i\,\hbox{\vrule height 8pt depth 8pt width 0.5pt}_{\,t=0}=
\nu(p)\cdot n^i(p)
\tag6.5
\endxalignat
$$
for Newtonian dynamical system \thetag{6.3}. Similar to initial data
\thetag{2.10} for modified Hamiltonian dynamical system, here initial
data \thetag{6.5} define {\it a shift} of hypersurface $\sigma$ along
trajectories of Newtonian dynamical system \thetag{6.3}. This shift
is called {\it normal shift} if hypersurfaces $\sigma_t$, which are
obtained from $\sigma$ by shift, keep orthogonality to shift trajectories
in time.
\definition{Definition 6.1} Newtonian dynamical system \thetag{6.3}
is called a system {\it admitting normal shift} if for any hypersurface
$\sigma$ there is a smooth function $\nu=\nu(p)$ on $\sigma$ such that
initial data \thetag{6.5} with this function $\nu$ define normal shift
of $\sigma$ along trajectories of dynamical system \thetag{6.3}.
\enddefinition
Suppose that $p_0$ is some fixed point of hypersurface $\sigma$ and
let $\nu_0$ be some fixed constant. Let's normalize $\nu(p)$ by
the following condition:
$$
\hskip -2em
\nu(p_0)=\nu_0.
\tag6.6
$$
\definition{Definition 6.2} Say that Newtonian dynamical system
\thetag{6.3} satisfies {\it strong normality condition} if for any
hypersurface $\sigma$, for any point $p_0\in\sigma$, and for
any constant $\nu_0\neq 0$ there is a smooth function $\nu=\nu(p)$
on $\sigma$ normalized by the condition \thetag{6.6} and such that
initial data \thetag{6.5} with this function $\nu$ define normal
shift of $\sigma$ along trajectories of dynamical system \thetag{6.3}.
\enddefinition
    Strong normality condition, in contrast to {\it the normality
condition} from definition~6.1, is less obvious. But it is more
convenient for to study by mathematical methods. In papers
\cite{12} and \cite{13} the following two systems of differential
equations for the force field $\bold F$ of Newtonian dynamical system
\thetag{6.3} were derived:
$$
\align
&\hskip -2em
\left\{
\aligned
&\sum^n_{i=1}\left(v^{-1}\,F_i+\shave{\sum^n_{j=1}}
\tilde\nabla_i\left(N^j\,F_j\right)\right)P^i_k=0,
\vspace{2ex}
&\aligned
&\sum^n_{i=1}\sum^n_{j=1}\left(\nabla_iF_j+\nabla_jF_i-2\,v^{-2}
\,F_i\,F_j\right)N^j\,P^i_k\,+\\
\vspace{0.4ex}
&+\sum^n_{i=1}\sum^n_{j=1}\left(\frac{F^j\,\tilde\nabla_jF_i}{v}
-\sum^n_{r=1}\frac{N^r\,N^j\,\tilde\nabla_jF_r}{v}\,F_i\right)P^i_k=0,
\endaligned
\endaligned\right.
\tag6.7\\
\vspace{2ex}
&\hskip -2em
\left\{
\aligned
&\aligned
&\sum^n_{i=1}\sum^n_{j=1}P^i_\varepsilon\,P^j_\sigma\left(
\,\shave{\sum^n_{m=1}}N^m\,\frac{F_i\,\tilde\nabla_mF_j}{v}
-\nabla_iF_j\right)=\\
\vspace{0.4ex}
&\quad\,=\sum^n_{i=1}\sum^n_{j=1}P^i_\varepsilon\,P^j_\sigma\left(
\,\shave{\sum^n_{m=1}}N^m\,\frac{F_j\,\tilde\nabla_mF_i}{v}
-\nabla_jF_i\right),
\endaligned\\
\vspace{1ex}
&\sum^n_{i=1}\sum^n_{j=1}P^j_\sigma\,\tilde\nabla_jF^i\,
P^\varepsilon_i=\sum^n_{i=1}\sum^n_{j=1}\sum^n_{m=1}
\frac{P^j_m\,\tilde\nabla_jF^i\,P^m_i}{n-1}\,P^\varepsilon_\sigma.
\endaligned\right.
\tag6.8
\endalign
$$
The equations \thetag{6.7} were called {\it weak normality equations},
while other equations \thetag{6.8} were called {\it additional normality
equations}. In Chapter~\uppercase\expandafter{\romannumeral 5} of thesis
\cite{23} the following theorem was proved.
\proclaim{Theorem 6.2} Newtonian dynamical system \thetag{6.3} satisfies
strong normality condition if and only if its force field $\bold F$
satisfies complete system of normality equations consisting of weak
normality equations \thetag{6.7} and including additional normality
equations \thetag{6.8} in the case of higher dimensions $n\geqslant 3$.
\endproclaim
    Weak normality equations \thetag{6.7} are related to {\it weak
normality condition}. In order to formulate this condition let's
consider one-parametric family of trajectories of Newtonian dynamical
system \thetag{6.3}. Denote it as follows:
$$
\hskip -2em
p=p(t,y).
\tag6.9
$$
Here $t$ is time variable and $y$ is a parameter. In local coordinates
this one-parametric family of trajectories \thetag{6.9} is expressed
by functions
$$
\hskip -2em
\cases
x^1=x^1(t,y),\\
.\ .\ .\ .\ .\ .\ .\ .\ .\ .\\
x^n=x^n(t,y).
\endcases
\tag6.10
$$
Differentiating \thetag{6.10} with respect to parameter $y$, we get
vector $\boldsymbol\tau$ with components
$$
\hskip -2em
\tau^i=\frac{\partial x^i}{\partial y}.
\tag6.11
$$
Vector $\boldsymbol\tau$ is called {\it vector of variation} of
trajectories. Then from \thetag{6.3} we derive
$$
\hskip -2em
\aligned
\nabla_{tt}\tau^k&=-\sum^n_{q=1}\sum^n_{i=1}\sum^n_{j=1}R^k_{qij}\,
\tau^i\,u^j\,u^q\,+\\
&+\sum^n_{q=1}\nabla_t\tau^q\,\tilde\nabla_qF^k+\sum^n_{q=1}\tau^q\,
\nabla_qF^k.
\endaligned
\tag6.12
$$
Here $R^k_{qij}$ components of curvature tensor for metric $\bold g$,
while $u^j$ and $u^q$ are components of velocity vector $\bold u$.
Function $\varphi$ defined as scalar product
$$
\hskip -2em
\varphi=(\bold u\,|\,\boldsymbol\tau)=\sum^n_{i=1}\sum^n_{j=1}
g_{ij}\,u^i\,\tau^j
\tag6.13
$$
is called {\it function of deviation}. From \thetag{6.12} one can
derive the following ordinary differential equation for the function
of deviation \thetag{6.13}:
$$
\hskip -2em
\sum^{2n}_{i=0}C_i(t)\,\varphi^{(i)}=0.
\tag6.14
$$
For general Newtonian dynamical system \thetag{6.3} this is
homogeneous ordinary differential equation of the order $2n$
(see details in Chapter~\uppercase\expandafter{\romannumeral 5}
of thesis \cite{23}). However, in special cases the equation
\thetag{6.14} can reduce to lower order differential equation.
Weak normality condition below specifies one of such cases.
Indeed, let's consider some trajectory $p=p(t)$ of Newtonian
dynamical system \thetag{6.3}. It can be included into
one-parametric family of trajectories \thetag{6.9} by various
ways. This defines various variation vectors $\boldsymbol
\tau$ with components satisfying differential equations
\thetag{6.12} and various deviation functions \thetag{6.13} on
the trajectory $p=p(t)$.
\definition{Definition 6.3} Say that Newtonian dynamical system
\thetag{6.3} satisfies {\it weak normality condition} if for each
its trajectory $p=p(t)$ and for any vector of variation $\boldsymbol
\tau$ on this trajectory corresponding function of deviation
$\varphi(t)$ satisfies homogeneous second order ordinary differential
equation 
$$
\hskip -2em
\ddot\varphi=\Cal A(t)\,\dot\varphi+\Cal B(t)\,\varphi
\tag6.15
$$
with coefficients depending only on choice of trajectory $p=p(t)$.
\enddefinition
\noindent As it was shown in paper \cite{12}, weak normality condition
is equivalent to weak normality equations \thetag{6.7} for the force
field of Newtonian dynamical system \thetag{6.3}.\par
    Now let's proceed with additional normality condition. In order to
formulate this condition let's consider some smooth hypersurface $\sigma$
in $M$ and let's fix some point $p_0$ on $\sigma$. Denote by $y^1,\,
\ldots,\,y^{n-1}$ local coordinates on $\sigma$ in some neighborhood of
fixed point $p_0$. Setting up initial data \thetag{6.5}, we can define
a family of trajectories $p=p(t,y^1,\ldots,y^n)$ of Newtonian dynamical
system \thetag{6.3} starting at the points of $\sigma$. Now this is
$(n-1)$-parametric family of trajectories expressed by functions
$$
\hskip -2em
\cases
x^1=x^1(t,y^1,\ldots,y^{n-1}),\\
.\ .\ .\ .\ .\ .\ .\ .\ .\ .\ .\ .\ .\ .\ .\ .\ .\ .\ .\\
x^n=x^n(t,y^1,\ldots,y^{n-1}).
\endcases
\tag6.16
$$
in local coordinates $x^1,\,\ldots,\,x^n$ in $M$. Differentiating
these functions \thetag{6.16} with respect to parameters $y^1,\,
\ldots,\,y^{n-1}$, as it is done in \thetag{6.11}, we get $n-1$
variation vectors $\boldsymbol\tau_1,\,\ldots,\,\boldsymbol\tau_{n-1}$.
It's easy to note that variation vectors $\boldsymbol\tau_1,\,\ldots,\,
\boldsymbol\tau_{n-1}$ form a frame of tangent vectors to hypersurfaces
$\sigma_t$ obtained by shifting initial hypersurface $\sigma$ along
trajectories \thetag{6.16} of dynamical system \thetag{6.3}. Therefore
corresponding deviation functions $\varphi_1,\,\ldots,\,\varphi_{n-1}$
serve as measure of orthogonality of $\sigma_t$ and shift trajectories
\thetag{6.16}. They should be identically zero in order to provide
orthogonality of shift: $\varphi_i(t,y^1,\ldots,y^{n-1})=0$. If we
consider initial data
$$
\xalignat 2
&\hskip -2em
\varphi_k\,\hbox{\vrule height 8pt depth 8pt width 0.5pt}_{\,t=0}=0,
&&\dot\varphi_k\,\hbox{\vrule height 8pt depth 8pt width 0.5pt}_{\,t=0}=0,
\tag6.17
\endxalignat
$$
then we can see that first part of initial conditions \thetag{16.17}
is fulfilled due to initial data \thetag{6.5}. Second part of these
conditions can be transformed to differential equations for the function
$\nu=\nu(p)=\nu(y^1,\ldots,y^{n-1})$ in \thetag{6.5}:
$$
\hskip -2em
\frac{\partial\nu}{\partial y^i}=-\nu^{-1}\,
(\bold F\,|\,\boldsymbol\tau_i).
\tag6.18
$$
If $\dim M=n\geqslant 3$, then the equations \thetag{6.18} form complete
system of Pfaff equations for scalar function $\nu$. The condition of
its compatibility is known as additional normality condition.
\definition{Definition 6.4} Say that Newtonian dynamical system
\thetag{6.3} satisfies {\it additional normality condition} if
for any smooth hypersurface $\sigma$ in $M$ and for any local
coordinates $y^1,\,\ldots,\,y^{n-1}$ on $\sigma$ corresponding
Pfaff equations \thetag{6.18} are compatible.
\enddefinition
   In paper \cite{13} it was shown that for $n\geqslant 3$ additional
normality condition is equivalent to additional normality equations
\thetag{6.8} for the force field of Newtonian dynamical system
\thetag{6.3}. In two-dimensional case $n=2$ situation is quite
different. Here we have only one parameter $y=y^1$ and \thetag{6.18}
turns to unique ordinary differential equation, which is compatible
with itself in anyway. Therefore in two-dimensional case additional
normality condition is always fulfilled. This special case is studied
in thesis \cite{24}.\par
    In higher dimensional case $n\geqslant 3$ complete system of
normality equations includes both \thetag{6.7} and \thetag{6.8}.
For this case in Chapter~\uppercase\expandafter{\romannumeral 7}
of thesis \cite{23} explicit formula for general solution of
complete system of normality equations was derived:
$$
\hskip -2em
F_k=\frac{h(W)\,N_k}{W'}
-|\bold u|\cdot\sum^n_{i=1}\frac{\nabla_iW}{W'}\cdot
\left(2\,N^i\,N_k-\delta^i_k\right).
\tag6.19
$$
Comparing formulas \thetag{6.4} and \thetag{6.19} we see that they
are quite similar. They differ only by first term in \thetag{6.19},
where $h=h(w)$ is arbitrary function of one variable. This fact
indicates that modified Lagrangian dynamical systems \thetag{5.2}
describing wave front dynamics and Newtonian dynamical systems
\thetag{6.3} admitting normal shift of hypersurfaces in Riemannian
manifolds are closely related with each other. In further sections
we are going to reveal this relation in more general case, when
manifold $M$ is not equipped with Riemannian metric. Problem of
interpreting first term in \thetag{6.19} should be considered in
separate paper.
\head
7. Weak normality phenomenon for modified Lagrangian
dynamical systems.
\endhead
    Let $M$ be a smooth manifold which is not equipped with
Riemannian metric, but which is equipped with modified Lagrangian
dynamical system. Let $L=L(p,\bold v)$ be Lagrange function for
this system. This is extended scalar field in contravariant
$\bold v$-representation. In the absence of Riemannian metric we
cannot use spatial gradient \thetag{3.4}. Therefore we write
modified Lagrange equations as
$$
\xalignat 2
&\hskip -2em
\dot x^i=\frac{v^i}{\Omega},
&&\frac{d}{dt}\left(\frac{\partial L}{\partial v^i}\right)=
\frac{1}{\Omega}\,\frac{\partial L}{\partial x^i}.
\tag7.1
\endxalignat
$$
Denominator $\Omega$ in \thetag{7.1} is determined by formula
\thetag{5.4}:
$$
\hskip -2em
\Omega=\sum^n_{i=1}v^i\,\tilde\nabla_iL=\sum^n_{i=1}v^i\,
\frac{\partial L}{\partial v^i}.
\tag7.2
$$
Remember that formula \thetag{4.2} defines Legendre map \thetag{4.3}.
Below we assume this map $\lambda$ to be invertible, and moreover,
we assume inverse map $\lambda^{-1}$ to be smooth. Local invertibility
of $\lambda$ means that matrix $\bold g$ with components
$$
\hskip -2em
\mu_{ij}=\frac{\tilde\nabla_{\!i}\tilde\nabla_{\!j}L}{2}=
\frac{1}{2}\,\frac{\partial^2L}{\partial v^i\,\partial v^j}
\tag7.3
$$
is non-degenerate. We shall assume this matrix to be positive:
$$
\hskip -2em
\boldsymbol\mu>0.
\tag7.4
$$
For real mechanical systems this condition $\boldsymbol\mu>0$ is fulfilled
since kinetic energy $K$ of such systems is positive quadratic function
of velocity vector. For such systems components of matrix $\boldsymbol\mu$
do not depend on $\bold v$, hence one can choose $\boldsymbol\mu$ to be
metric tensor for Riemannian metric in $M$. However, we shall consider
more general case, when $\boldsymbol\mu$ is extended tensor field with
components depending on $\bold v$.\par
   In addition to inequality \thetag{7.4} we shall assume that
denominator $\Omega$ in modified Lagrange equations \thetag{7.1}
(which is determined by \thetag{7.2}) is positive function:
$$
\hskip -2em
\Omega>0\text{, \ for \ }\bold v\neq 0.
\tag7.5
$$
This assumption is consistent since for real mechanical systems $\Omega=
2\,K$.\par
   Now let's consider one-parametric family of trajectories
$p=p(t,y)$ of modified Lagrangian dynamical system \thetag{7.1}.
In local coordinates these curves are expressed by functions
\thetag{6.10}. Formula \thetag{6.11} then defines vector of
variation $\boldsymbol\tau$. In order to define function of
deviation $\varphi$ we could use formula \thetag{6.13} with matrix
\thetag{7.3} as metric. However, we choose another formula for
$\varphi$:
$$
\hskip -2em
\varphi=\left<\bold p\,|\,\boldsymbol\tau\right>=
\sum^n_{i=1}\tilde\nabla_{\!i}L\ \tau^i.
\tag7.6
$$
Here $\bold p$ is momentum covector defined by formula \thetag{4.2},
while angular brackets denote contraction of vector and covector
\footnotemark.\footnotetext{Such notations are often used in quantum
mechanics. See \cite{31}.}\adjustfootnotemark{-1} Function $\varphi$
in \thetag{7.1} can
be treated as scalar product of vectors $\bold v$ and $\boldsymbol\tau$.
This scalar product is linear with respect to vector $\boldsymbol
\tau$, but it is nonlinear respect to vector $\bold v$. Such scalar
products usually arise in Finslerian geometry (see
Chapter~\uppercase\expandafter{\romannumeral 8} of thesis \cite{23}).
\par
   Let $\bold v=\bold v(t,y)$ be velocity vector for one-parametric
family of trajectories $p=p(t,y)$ of modified Lagrangian dynamical
system \thetag{7.1}. In local coordinates this vector-function is
expressed by the following scalar functions:
$$
\hskip -2em
\cases
v^1=v^1(t,y),\\
.\ .\ .\ .\ .\ .\ .\ .\ .\ .\\
v^n=v^n(t,y).
\endcases
\tag7.7
$$
Differentiating \thetag{7.7} with respect to parameter $y$, we get
series of functions
$$
\hskip -2em
\theta^i=\frac{\partial v^i}{\partial y}=\dot\tau^i.
\tag7.8
$$
In contrast to $\tau^i$ in \thetag{6.11}, these functions $\theta^1,\,
\ldots,\,\theta^n$ are not interpreted as components of vector. We
shall use them in order to simplify further calculations. Differentiating
first equation \thetag{7.1} with respect to $y$, we obtain
$$
\hskip -2em
\gathered
\dot\tau^i=\frac{\theta^i}{\Omega}-\sum^n_{s=1}\frac{v^i}
{\Omega^2}\,\frac{\partial L}{\partial v^s}\,\theta^s\,-\\
\vspace{1ex}
-\sum^n_{k=1}\sum^n_{s=1}\frac{\partial^2\!L}{\partial v^k\,
\partial v^s}\,\frac{v^i\,v^k\,\theta^s}{\Omega^2}
-\sum^n_{k=1}\sum^n_{s=1}\frac{\partial^2\!L}{\partial v^k\,
\partial x^s}\,\frac{v^i\,v^k\,\tau^s}{\Omega^2}.
\endgathered
\tag7.9
$$
Differentiating second equation \thetag{7.1} with respect to
parameter $y$, we get
$$
\hskip -2em
\gathered
\sum^n_{s=1}\frac{\partial^2\!L}{\partial v^i\,\partial v^s}\,
\dot\theta^s+\sum^n_{s=1}\frac{\partial^2\!L}{\partial v^i\,
\partial x^s}\,\dot\tau^s+
\sum^n_{k=1}\sum^n_{s=1}\frac{\partial^3\!L}{\partial v^i\,
\partial v^s\,\partial v^k}\,\dot v^k\,\theta^s\,+\\
+\sum^n_{k=1}\sum^n_{s=1}\frac{\partial^3\!L}{\partial v^i\,
\partial v^s\,\partial x^k}\,\dot x^k\,\theta^s
+\sum^n_{k=1}\sum^n_{s=1}\frac{\partial^3\!L}{\partial v^i\,
\partial x^s\,\partial v^k}\,\dot v^k\,\tau^s\,+\\
\vspace{1ex}
+\sum^n_{k=1}\sum^n_{s=1}\frac{\partial^3\!L}{\partial v^i\,
\partial x^s\,\partial x^k}\,\dot x^k\,\tau^s
=\sum^n_{s=1}\frac{1}{\Omega}\,\frac{\partial^2\!L}{\partial x^i\,
\partial v^s}\,\theta^s\,+\\
\vspace{1ex}
+\sum^n_{s=1}\frac{1}{\Omega}\,\frac{\partial^2\!L}{\partial x^i\,
\partial x^s}\,\tau^s-\sum^n_{s=1}\frac{\partial L}{\partial x^i}\,
\frac{\partial L}{\partial v^s}\,\frac{\theta^s}{\Omega^2}\,-\\
\vspace{1ex}
-\sum^n_{k=1}\sum^n_{s=1}\frac{\partial L}{\partial x^i}\,
\frac{\partial^2\!L}{\partial v^k\,\partial v^s}
\frac{\,v^k\,\theta^s}{\Omega^2}
-\sum^n_{k=1}\sum^n_{s=1}\frac{\partial L}{\partial x^i}\,
\frac{\partial^2\!L}{\partial v^k\,\partial x^s}
\frac{\,v^k\,\tau^s}{\Omega^2}.
\endgathered
\tag7.10
$$
Both \thetag{7.9} and \thetag{7.10} form a system of homogeneous
linear ordinary differential equations with respect to functions
$\tau^1,\,\ldots,\,\tau^n$ and $\theta^1,\,\ldots,\,\theta^n$.
This system of equations is an analog of equations \thetag{6.12}
considered above.\par
    Function of deviation $\varphi$ defined by formula \thetag{7.6}
depends linearly on components of vector $\boldsymbol\tau$. Let's
calculate time derivatives of this function. For $\dot\varphi$ we get
$$
\gathered
\dot\varphi=\sum^n_{s=1}\frac{\partial L}{\partial v^s}\,\dot\tau^s
+\sum^n_{k=1}\frac{d}{dt}\!\left(\frac{\partial L}{\partial v^s}\right)
\tau^s=\sum^n_{s=1}\frac{\partial L}{\partial v^s}\,\dot\tau^s\,+\\
\vspace{1ex}
+\sum^n_{s=1}\frac{\partial L}{\partial x^s}\,\frac{\tau^s}{\Omega}
=\sum^n_{s=1}\frac{\partial L}{\partial v^s}\,\frac{\theta^s}{\Omega}-
\sum^n_{s=1}\sum^n_{r=1}\frac{\partial L}{\partial v^s}\,\frac{v^s}
{\Omega^2}\,\frac{\partial L}{\partial v^r}\,\theta^r\,-\\
\vspace{1ex}
-\sum^n_{s=1}\sum^n_{k=1}\sum^n_{r=1}\frac{\partial L}{\partial v^s}\,
\frac{\partial^2\!L}{\partial v^k\,
\partial v^r}\,\frac{v^s\,v^k\,\theta^r}{\Omega^2}
-\sum^n_{s=1}\sum^n_{k=1}\sum^n_{r=1}\frac{\partial L}{\partial v^s}\,
\frac{\partial^2\!L}{\partial v^k\,
\partial x^r}\,\frac{v^s\,v^k\,\tau^r}{\Omega^2}\,+\\
\vspace{1ex}
+\sum^n_{s=1}\frac{\partial L}{\partial x^s}\,\frac{\tau^s}{\Omega}=
\sum^n_{s=1}\frac{\partial L}{\partial x^s}\,\frac{\tau^s}{\Omega}
-\sum^n_{k=1}\sum^n_{r=1}\frac{\partial^2\!L}{\partial v^k\,
\partial v^r}\,\frac{v^k\,\theta^r}{\Omega}
-\sum^n_{k=1}\sum^n_{r=1}\frac{\partial^2\!L}{\partial v^k\,
\partial x^r}\,\frac{v^k\,\tau^r}{\Omega}.
\endgathered
$$
In the above calculations we used second equation \thetag{7.1},
used formula \thetag{7.2} for denominator $\Omega$, and used
differential equations \thetag{7.9}. As a result we obtained
$$
\dot\varphi=\sum^n_{s=1}\frac{\partial L}{\partial x^s}\,
\frac{\tau^s}{\Omega}
-\sum^n_{k=1}\sum^n_{r=1}\frac{\partial^2\!L}{\partial v^k\,
\partial v^r}\,\frac{v^k\,\theta^r}{\Omega}
-\sum^n_{k=1}\sum^n_{r=1}\frac{\partial^2\!L}{\partial v^k\,
\partial x^r}\,\frac{v^k\,\tau^r}{\Omega}.
\tag7.11
$$
Now let's differentiate  \thetag{7.11} once more. This yields
$$
\gather
\ddot\varphi=
-\sum^n_{i=1}\sum^n_{s=1}\frac{\partial^2\!L}{\partial v^i\,
\partial v^s}\,\frac{v^i\,\dot\theta^s}{\Omega}
-\sum^n_{i=1}\sum^n_{s=1}\frac{\partial^2\!L}{\partial v^i\,
\partial x^s}\,\frac{v^i\,\dot\tau^s}{\Omega}
+\sum^n_{s=1}\frac{\partial L}{\partial x^s}\,\frac{\dot\tau^s}
{\Omega}-\\
\vspace{1ex}
-\sum^n_{k=1}\sum^n_{i=1}\sum^n_{s=1}\frac{\partial^3\!L}
{\partial v^k\,\partial v^i\,\partial v^s}\,\frac{\dot v^k\,v^i\,
\theta^s}{\Omega}
-\sum^n_{k=1}\sum^n_{i=1}\sum^n_{s=1}\frac{\partial^2\!L}
{\partial x^k\,\partial v^i\,\partial v^s}\,\frac{\dot x^k\,v^i\,
\theta^s}{\Omega}\,-\\
\vspace{1ex}
-\sum^n_{k=1}\sum^n_{i=1}\sum^n_{s=1}\frac{\partial^3\!L}
{\partial v^k\,\partial v^i\,\partial x^s}\,\frac{\dot v^k\,v^i\,
\tau^s}{\Omega}
-\sum^n_{k=1}\sum^n_{i=1}\sum^n_{s=1}\frac{\partial^3\!L}
{\partial x^k\,\partial v^i\,\partial x^s}\,\frac{\dot x^k\,v^i\,
\tau^s}{\Omega}\,-\\
\vspace{1ex}
-\sum^n_{i=1}\sum^n_{s=1}\frac{\partial^2\!L}{\partial v^i\,
\partial v^s}\,\frac{\dot v^i\,\theta^s}{\Omega}
-\sum^n_{i=1}\sum^n_{s=1}\frac{\partial^2\!L}{\partial v^i\,
\partial x^s}\,\frac{\dot v^i\,\tau^s}{\Omega}\,+\\
\vspace{1ex}
+\sum^n_{k=1}\sum^n_{s=1}\frac{\partial^2\!L}{\partial v^k\,
\partial x^s}\,\frac{\dot v^k\,\tau^s}{\Omega}
+\sum^n_{k=1}\sum^n_{s=1}\frac{\partial^2\!L}{\partial x^k\,
\partial x^s}\,\frac{\dot x^k\,\tau^s}{\Omega}-\frac{\dot\Omega}
{\Omega}\,\dot\varphi.
\endgather
$$
Using equations \thetag{7.10}, we can eliminate all entries of
derivatives $\dot\theta^s$ from the above expression for $\ddot\varphi$.
As a result we get reduced formula
$$
\gather
\ddot\varphi=
-\sum^n_{i=1}\sum^n_{s=1}\frac{\partial^2\!L}{\partial x^i\,
\partial v^s}\,\frac{v^i\,\theta^s}{\Omega^2}
-\sum^n_{i=1}\sum^n_{s=1}\frac{\partial^2\!L}{\partial x^i\,
\partial x^s}\,\frac{v^i\,\tau^s}{\Omega^2}\,+\\
\displaybreak
+\sum^n_{i=1}\sum^n_{s=1}\frac{\partial L}{\partial x^i}\,
\frac{\partial L}{\partial v^s}\,\frac{v^i\,\theta^s}{\Omega^3}
+\sum^n_{i=1}\sum^n_{k=1}\sum^n_{s=1}\frac{\partial L}
{\partial x^i}\,\frac{\partial^2\!L}{\partial v^k\,\partial v^s}
\frac{v^i\,v^k\,\theta^s}{\Omega^3}\,+\\
\vspace{1ex}
+\sum^n_{i=1}\sum^n_{k=1}\sum^n_{s=1}\frac{\partial L}
{\partial x^i}\,\frac{\partial^2\!L}{\partial v^k\,\partial x^s}
\frac{v^i\,v^k\,\tau^s}{\Omega^3}
+\sum^n_{s=1}\frac{\partial L}{\partial x^s}\,\frac{\dot\tau^s}
{\Omega}\,-\\
\vspace{1ex}
-\sum^n_{i=1}\sum^n_{s=1}\frac{\partial^2\!L}{\partial v^i\,
\partial v^s}\,\frac{\dot v^i\,\theta^s}{\Omega}
+\sum^n_{k=1}\sum^n_{s=1}\frac{\partial^2\!L}{\partial x^k\,
\partial x^s}\,\frac{\dot x^k\,\tau^s}{\Omega}-\frac{\dot\Omega}
{\Omega}\,\dot\varphi.
\endgather
$$
Now we eliminate all entries of $\dot\tau^s$, using equations
\thetag{7.9} for this purpose:
$$
\hskip +0.5em
\gathered
\ddot\varphi=
-\sum^n_{i=1}\sum^n_{s=1}\frac{\partial^2\!L}{\partial x^i\,
\partial v^s}\,\frac{v^i\,\theta^s}{\Omega^2}
-\sum^n_{i=1}\sum^n_{s=1}\frac{\partial^2\!L}{\partial x^i\,
\partial x^s}\,\frac{v^i\,\tau^s}{\Omega^2}\,+\\
\vspace{1ex}
+\sum^n_{s=1}\frac{\partial L}{\partial x^s}\,\frac{\theta^s}
{\Omega^2}-\sum^n_{i=1}\sum^n_{s=1}\frac{\partial^2\!L}
{\partial v^i\,\partial v^s}\,\frac{\dot v^i\,\theta^s}{\Omega}
+\sum^n_{k=1}\sum^n_{s=1}\frac{\partial^2\!L}{\partial x^k\,
\partial x^s}\,\frac{\dot x^k\,\tau^s}{\Omega}-\frac{\dot\Omega}
{\Omega}\,\dot\varphi.
\endgathered
\tag7.12
$$
Then, using first equation in \thetag{7.1}, we express time derivative
$\dot x^k$ through $v^k$. As a result formula \thetag{7.12}
for $\ddot\varphi$ reduces to the following one:
$$
\hskip +0.5em
\gathered
\ddot\varphi=
-\sum^n_{i=1}\sum^n_{s=1}\frac{\partial^2\!L}{\partial x^i\,
\partial v^s}\,\frac{v^i\,\theta^s}{\Omega^2}
+\sum^n_{s=1}\frac{\partial L}{\partial x^s}\,\frac{\theta^s}
{\Omega^2}\,-\\
-\sum^n_{i=1}\sum^n_{s=1}\frac{\partial^2\!L}{\partial v^i\,
\partial v^s}\,\frac{\dot v^i\,\theta^s}{\Omega}
-\frac{\dot\Omega}{\Omega}\,\dot\varphi.
\endgathered
\tag7.13
$$
Now, if one take into account second equation \thetag{7.1} written
in expanded form, then three terms in right hand side of
\thetag{7.13} can be canceled. This yields
$$
\hskip -2em
\ddot\varphi+\frac{\dot\Omega}{\Omega}\,\dot\varphi=0.
\tag7.14
$$
Thus, in the end of huge calculations we get very simple relationship
\thetag{7.14}, which is homogeneous second order linear ordinary
differential equation. It is even more simple than analogous equation
\thetag{6.15} considered in previous section. Now, if we formulate
definition~6.3 respective to modified Lagrangian dynamical system
\thetag{7.1} and if we use formula \thetag{7.6} for $\varphi$, then
from \thetag{7.14} we derive the following theorem.
\proclaim{Theorem 7.1} Each modified Lagrangian dynamical system
\thetag{7.1} satisfies {\bf weak normality condition} with respect
to deviation functions \thetag{7.6} determined by its own Lagrange
function $L$.
\endproclaim
\head
8. Additional normality phenomenon.
\endhead
    In order to reproduce results of section~6 in present more
complicated geometric environment we should consider some hypersurface
$\sigma$ in $M$, and we should arrange a shift of $\sigma$ by means
of modified Lagrangian dynamical system \thetag{7.1}. Fortunately we
should not invent something absolutely new for this purpose. Wave
front dynamics considered in section~2 suggests a way of how to do
this. In the absence of Riemannian metric we cannot choose unit normal
vector on $\sigma$. However, we can take normal covector $\bold n=
\bold n(p)$, which is unique up to a scalar factor. Then we can set up
Cauchy problem with the following initial data:
$$
\xalignat 2
&\hskip -2em
x^i\,\hbox{\vrule height 8pt depth 8pt width 0.5pt}_{\,t=0}
=x^i(p),
&&p_i\,\hbox{\vrule height 8pt depth 8pt width 0.5pt}_{\,t=0}=
\nu(p)\cdot n_i(p)
\tag8.1
\endxalignat
$$
(compare with \thetag{2.10} above). Here $p$ is a point of $\sigma$ and
$p_i$ are components of momentum covector $\bold p$ defined by formula
\thetag{4.2}:
$$
p_i=\tilde\nabla_iL=\frac{\partial L}{\partial v^i}.
$$
Initial data \thetag{8.1} determine initial velocity $\bold v$
implicitly through initial momentum covector $\bold p$ due to
invertibility of Legendre map $\lambda$. Applying initial data
\thetag{8.1} to modified Lagrangian dynamical system \thetag{7.1},
we obtain a family of trajectories $p=p(t,y^1,\ldots,y^n)$ starting
at the points of hypersurface $\sigma$. Similar to \thetag{6.16},
in local coordinates these trajectories are expressed by the following
functions:
$$
\hskip -2em
\cases
x^1=x^1(t,y^1,\ldots,y^{n-1}),\\
.\ .\ .\ .\ .\ .\ .\ .\ .\ .\ .\ .\ .\ .\ .\ .\ .\ .\\
x^n=x^n(t,y^1,\ldots,y^{n-1}).
\endcases
\tag8.2
$$
These functions \thetag{8.2} determine variation vectors $\boldsymbol
\tau_1,\,\ldots,\,\boldsymbol\tau_{n-1}$ with components
$$
\tau^i_j=\frac{\partial x^i}{\partial y^j}
$$
(compare with \thetag{6.11} and \thetag{7.8}). Each variation vector
determines corresponding deviation function according to the formula
\thetag{7.6}. We denote these deviation functions by $\varphi_1,\,
\ldots,\,\varphi_{n-1}$ as in section~6.
\definition{Definition 8.1} Shift of initial hypersurface $\sigma$
determined by modified Lagrangian dynamical system \thetag{7.1} and
by initial data \thetag{8.1} for it is called {\it normal shift} in
inner geometry of dynamical system \thetag{7.1} if all deviation
functions $\varphi_1,\,\ldots,\,\varphi_{n-1}$ are identically zero.
\enddefinition
    Due to differential equation \thetag{7.14} for deviation functions
in order to arrange a normal shift of $\sigma$ it is sufficient to
provide initial conditions
$$
\xalignat 2
&\hskip -2em
\varphi_i,\hbox{\vrule height 8pt depth 8pt width 0.5pt}_{\,t=0}=0,
&&\dot\varphi_i\,\hbox{\vrule height 8pt depth 8pt width 0.5pt}_{\,t=0}=0
\tag8.3
\endxalignat
$$
just the same as in \thetag{6.17}. First part of initial conditions
\thetag{8.3} is fulfilled due to initial data \thetag{8.1}. Second
part of these conditions should be transformed to differential equations
for the function $\nu=\nu(p)=\nu(y^1,\ldots,y^{n-1})$ in \thetag{8.1}.
For this purpose we could use formula \thetag{7.11} derived in section~7.
However, initial data \thetag{8.1} explicitly relate function $\nu=\nu(p)$
with initial value of momentum covector $\bold p$, while relation to
velocity vector $\bold v$ is implicit. Therefore it is easier to transform
formula \thetag{7.6} to $\bold p$-representation. This yields
$$
\hskip -2em
\varphi_i=\left<\bold p\,|\,\boldsymbol\tau_i\right>=
\sum^n_{s=1}p_s\,\tau^s_i.
\tag8.4
$$
Remember that modified Lagrange equations, when transformed to
$\bold p$-representa\-tion, look like modified Hamilton equations
\thetag{2.8}:
$$
\xalignat 2
&\hskip -2em
\dot x^i=\frac{1}{\Omega}\,\frac{\partial H}{\partial p_i},
&&\dot p_i=-\frac{1}{\Omega}\,\frac{\partial H}{\partial x^i}.
\tag8.5
\endxalignat
$$
Here Hamilton function $H$ is determined by formula \thetag{4.5}
and $\Omega$ is given by formula \thetag{5.3}. Now, differentiating
formula \thetag{8.4}, we obtain
$$
\dot\varphi_i=\sum^n_{s=1}p_s\,\tau^s_i=\sum^n_{s=1}\dot p_s\,\tau^s_i
+\sum^n_{s=1}p_s\,\dot\tau^s_i=-\sum^n_{s=1}\frac{1}{\Omega}\,
\frac{\partial H}{\partial x^s}\,\tau^s_i+\sum^n_{s=1}p_s\,\dot\tau^s_i.
\tag8.6
$$
In order to calculate time derivatives $\dot\tau^s_i$ in formula
\thetag{8.6} we use first part of modified Hamilton equations
\thetag{8.5}. A a result for $\dot\tau^s_i$ we get
$$
\gather
\dot\tau^s_i=\frac{\partial^2 x^s}{\partial t\,\partial y^i}=
\frac{\partial}{\partial y^i}\!\left(\frac{1}{\Omega}\,\frac{\partial H}
{\partial p_s}\right)=
\sum^n_{r=1}\frac{\partial x^r}{\partial y^i}\,
\frac{\partial}{\partial x^r}\!\left(\frac{1}{\Omega}\,\frac{\partial H}
{\partial p_s}\right)\,+\\
\vspace{1ex}
+\sum^n_{r=1}\frac{\partial p_r}{\partial y^i}\,
\frac{\partial}{\partial p_r}\!\left(\frac{1}{\Omega}\,\frac{\partial H}
{\partial p_s}\right)=
\sum^n_{r=1}\tau^r_i\,\frac{\partial}{\partial x^r}\!\left(\frac{1}
{\Omega}\,\frac{\partial H}{\partial p_s}\right)
+\sum^n_{r=1}\frac{\partial p_r}{\partial y^i}\,
\frac{\partial}{\partial p_r}\!\left(\frac{1}{\Omega}\,\frac{\partial H}
{\partial p_s}\right).
\endgather
$$
Let's substitute this expression into \thetag{8.6}. This yields
$$
\gather
\dot\varphi_i=-\sum^n_{s=1}\frac{1}{\Omega}\,
\frac{\partial H}{\partial x^s}\,\tau^s_i+\sum^n_{s=1}\sum^n_{r=1}
p_s\,\tau^r_i\,\frac{\partial}{\partial x^r}\!\left(\frac{1}
{\Omega}\,\frac{\partial H}{\partial p_s}\right)\,+\\
\vspace{1ex}
+\sum^n_{s=1}\sum^n_{r=1}p_s\,\frac{\partial p_r}{\partial y^i}\,
\frac{\partial}{\partial p_r}\!\left(\frac{1}{\Omega}\,\frac{\partial H}
{\partial p_s}\right)=\sum^n_{r=1}
\tau^r_i\,\frac{\partial}{\partial x^r}\!\left(\,\shave{\sum^n_{s=1}}
\frac{p_s}{\Omega}\,\frac{\partial H}{\partial p_s}\right)\,+\\
\vspace{1ex}
+\sum^n_{r=1}\frac{\partial p_r}{\partial y^i}\,
\frac{\partial}{\partial p_r}\!\left(\,\shave{\sum^n_{s=1}}\frac{p_s}
{\Omega}\,\frac{\partial H}{\partial p_s}\right)
-\sum^n_{s=1}\frac{1}{\Omega}\,\frac{\partial H}{\partial p_s}\,
\frac{\partial p_s}{\partial y^i}-\sum^n_{s=1}\frac{1}{\Omega}\,
\frac{\partial H}{\partial x^s}\,\tau^s_i.
\endgather
$$
First two terms in right hand side of the above equality are identically
zero. This follows from formula \thetag{5.3} for $\Omega$. Thus for
$\dot\varphi_i$ we have
$$
\hskip -2em
\dot\varphi_i=-\sum^n_{s=1}\frac{1}{\Omega}\,\frac{\partial H}
{\partial p_s}\,\frac{\partial p_s}{\partial y^i}
-\sum^n_{s=1}\frac{1}{\Omega}\,\frac{\partial H}{\partial x^s}
\,\tau^s_i.
\tag8.7
$$
If we recall initial conditions \thetag{8.3}, then from \thetag{8.7}
we derive
$$
\hskip -2em
\sum^n_{s=1}\frac{\partial H}{\partial x^s}
\,\tau^s_i+
\sum^n_{s=1}\frac{\partial H}{\partial p_s}\,
\left(\frac{\partial p_s}{\partial y^i}\right)\!
\,\hbox{\vrule height 12pt depth 8pt width 0.5pt}_{\,t=0}=0.
\tag8.8
$$
Calculating partial derivatives $\partial p_s/\partial y^i$ in
\thetag{8.8}, we should remember \thetag{8.1}. Then
$$
\hskip -2em
\left(\frac{\partial p_s}{\partial y^i}\right)\!
\,\hbox{\vrule height 12pt depth 8pt width 0.5pt}_{\,t=0}
=\frac{\partial\nu}{\partial y^i}\,n_s+
\nu\,\frac{\partial n_s}{\partial y^i}
=\frac{1}{\nu}\,\frac{\partial\nu}{\partial y^i}\,p_s+
\nu\,\frac{\partial n_s}{\partial y^i}.
\tag8.9
$$
Substituting this expression into \thetag{8.8} and using formula
\thetag{5.3} for $\Omega$, we can transform \thetag{8.8} to the
partial differential equations for $\nu$:
$$
\hskip -2em
\frac{1}{\nu}\,\frac{\partial\nu}{\partial y^i}=
-\sum^n_{s=1}\frac{\nu}{\Omega}\,\frac{\partial n_s}{\partial y^i}\,
\frac{\partial H}{\partial p_s}
-\sum^n_{s=1}\frac{\partial H}{\partial x^s}\,\frac{\tau^s_i}{\Omega}.
\tag8.10
$$
Differential equations \thetag{8.10} are analogs of the equations
\thetag{6.18}. If $n\geqslant 3$, then they form complete
system of Pfaff equations for the function $\nu=
\nu(y^1,\ldots,y^{n-1})$. Therefore we can formulate {\it additional
normality condition} for modified Lagrangian dynamical system \thetag{7.1}
as compatibility condition for Pfaff equations \thetag{8.10}.\par
    Suppose that $n\geqslant 3$. Let's examine if differential equations
\thetag{8.10} are compatible. For this purpose let's calculate second
order partial derivatives of $\nu$ using \thetag{8.10}:
$$
\gather
\frac{\partial^2\nu}{\partial y^i\,\partial y^j}=-\sum^n_{s=1}
\frac{\nu^2}{\Omega}\,\frac{\partial H}{\partial p_s}\,
\frac{\partial^2 n_s}{\partial y^i\,\partial y^j}-\sum^n_{s=1}
\frac{\partial H}{\partial x^s}\,\frac{\nu}{\Omega}\,\frac{\partial^2 x^s}
{\partial y^i\,\partial y^j}\,+\\
\vspace{1ex}
+\,\frac{2\,\nu^3}{\Omega^2}\sum^n_{s=1}\sum^n_{r=1}\frac{\partial H}
{\partial p_s}\,\frac{\partial H}{\partial p_r}\,\frac{\partial n_s}
{\partial y^i}\,\frac{\partial n_r}{\partial y^j}
+\frac{\nu}{\Omega^2}\sum^n_{s=1}\sum^n_{r=1}\frac{\partial H}
{\partial x^s}\,\frac{\partial H}{\partial x^r}\,\tau^s_i\,\tau^r_j\,+\\
\vspace{1ex}
+\,\frac{2\,\nu^2}{\Omega^2}\sum^n_{s=1}\sum^n_{r=1}\frac{\partial H}
{\partial x^s}\,\frac{\partial H}{\partial p_r}\,\frac{\partial n_r}
{\partial y^j}\,\tau^s_i
+\frac{\nu^2}{\Omega^2}\sum^n_{s=1}\sum^n_{r=1}\frac{\partial H}
{\partial p_s}\,\frac{\partial H}{\partial x^r}\,
\frac{\partial n_s}{\partial y^i}\,\tau^r_j\,-\\
\vspace{1ex}
-\sum^n_{s=1}\sum^n_{r=1}\nu^2\,\frac{\partial}{\partial p_s}\!
\left(\frac{1}{\Omega}\,\frac{\partial H}{\partial p_r}\right)
\frac{\partial p_s}{\partial y^i}\,\frac{\partial n_r}{\partial y^j}
-\sum^n_{s=1}\sum^n_{r=1}\nu\,\,\frac{\partial}{\partial p_s}\!
\left(\frac{1}{\Omega}\,\frac{\partial H}{\partial x^r}\right)
\frac{\partial p_s}{\partial y^i}\,\tau^r_j\,-\\
\vspace{1ex}
-\sum^n_{s=1}\sum^n_{r=1}\nu^2\,\frac{\partial}{\partial x^s}\!
\left(\frac{1}{\Omega}\,\frac{\partial H}{\partial p_r}\right)
\frac{\partial n_r}{\partial y^j}\,\tau^s_i
-\sum^n_{s=1}\sum^n_{r=1}\nu\,\,\frac{\partial}{\partial x^s}\!
\left(\frac{1}{\Omega}\,\frac{\partial H}{\partial x^r}\right)
\tau^s_i\,\tau^r_j.
\endgather
$$
For Pfaff equations \thetag{8.10} to be compatible, right hand side
of the above equality should be symmetric in indices $i$ and $j$.
First four terms there are obviously symmetric. Below we shall not
write such terms explicitly denoting them by dots:
$$
\gather
\frac{\partial^2\nu}{\partial y^i\,\partial y^j}=\dots
+\sum^n_{s=1}\sum^n_{r=1}\left(\frac{2\,\nu^2}{\Omega^2}\,
\frac{\partial H}{\partial x^s}\,\frac{\partial H}{\partial p_r}
-\nu^2\,\frac{\partial}{\partial x^s}\!
\left(\frac{1}{\Omega}\,\frac{\partial H}{\partial p_r}\right)
\right)\frac{\partial n_r}{\partial y^j}\,\tau^s_i\,+\\
\vspace{1ex}
+\sum^n_{s=1}\sum^n_{r=1}\left(\frac{\nu^2}{\Omega^2}\,
\frac{\partial H}{\partial p_s}\,\frac{\partial H}{\partial x^r}-
\nu^2\,\frac{\partial}{\partial p_s}\!
\left(\frac{1}{\Omega}\,\frac{\partial H}{\partial x^r}\right)
\right)
\frac{\partial n_s}{\partial y^i}\,\tau^r_j\,-\\
\vspace{1ex}
-\sum^n_{s=1}\sum^n_{r=1}\nu\,p_s\,\frac{\partial}{\partial p_s}\!
\left(\frac{1}{\Omega}\,\frac{\partial H}{\partial p_r}\right)
\frac{\partial\nu}{\partial y^i}\,\frac{\partial n_r}{\partial y^j}
-\sum^n_{s=1}\sum^n_{r=1}p_s\,\frac{\partial}{\partial p_s}\!
\left(\frac{1}{\Omega}\,\frac{\partial H}{\partial x^r}\right)
\frac{\partial\nu}{\partial y^i}\,\tau^r_j\,-\\
\vspace{1ex}
-\sum^n_{s=1}\sum^n_{r=1}\nu^3\,\frac{\partial}{\partial p_s}\!
\left(\frac{1}{\Omega}\,\frac{\partial H}{\partial p_r}\right)
\frac{\partial n_s}{\partial y^i}\,\frac{\partial n_r}{\partial y^j}
-\sum^n_{s=1}\sum^n_{r=1}\nu\,\,\frac{\partial}{\partial x^s}\!
\left(\frac{1}{\Omega}\,\frac{\partial H}{\partial x^r}\right)
\tau^s_i\,\tau^r_j.
\endgather
$$
In the above calculations we used formula \thetag{8.9} for partial
derivatives $\partial p_s/\partial y^i$. Below we use the equations
\thetag{8.10} for to express partial derivatives $\partial\nu/
\partial y^i$:
$$
\gather
\frac{\partial^2\nu}{\partial y^i\,\partial y^j}=\dots
+\sum^n_{s=1}\sum^n_{r=1}\left(\frac{2\,\nu^2}{\Omega^2}\,
\frac{\partial H}{\partial x^s}\,\frac{\partial H}{\partial p_r}
-\nu^2\,\frac{\partial}{\partial x^s}\!
\left(\frac{1}{\Omega}\,\frac{\partial H}{\partial p_r}\right)\,+
\vphantom{\shave{\sum^n_{q=1}}}\right.\\
\vspace{1ex}
\left.+\shave{\sum^n_{q=1}}\frac{\nu^2\,p_q}{\Omega}\,\frac{\partial}
{\partial p_q}\!\left(\frac{1}{\Omega}\,\frac{\partial H}
{\partial p_r}\right)\frac{\partial H}{\partial x^s}
\right)\frac{\partial n_r}{\partial y^j}\,\tau^s_i
+\sum^n_{s=1}\sum^n_{r=1}\left(\frac{\nu^2}{\Omega^2}\,
\frac{\partial H}{\partial p_s}\,\frac{\partial H}{\partial x^r}\,-
\vphantom{\shave{\sum^n_{q=1}}}\right.\\
\vspace{1ex}
\left.-\nu^2\,\frac{\partial}{\partial p_s}\!
\left(\frac{1}{\Omega}\,\frac{\partial H}{\partial x^r}\right)
+\shave{\sum^n_{q=1}}\frac{\nu^2\,p_q}{\Omega}\,
\frac{\partial}{\partial p_q}\!\left(\frac{1}{\Omega}\,
\frac{\partial H}{\partial x^r}\right)\frac{\partial H}
{\partial p_s}\right)
\frac{\partial n_s}{\partial y^i}\,\tau^r_j\,+\\
\vspace{1ex}
+\sum^n_{s=1}\sum^n_{r=1}\left(\,\shave{\sum^n_{q=1}}
\frac{\nu^3\,p_q}{\Omega}\,\frac{\partial}{\partial p_q}\!
\left(\frac{1}{\Omega}\,\frac{\partial H}{\partial p_r}\right)
\frac{\partial H}{\partial p_s}+\frac{\nu^3}{\Omega^2}\,  
\frac{\partial\Omega}{\partial p_s}\,\frac{\partial H}
{\partial p_r}\right)\frac{\partial n_s}{\partial y^i}\,
\frac{\partial n_r}{\partial y^j}\,+\\
\vspace{1ex}
+\sum^n_{s=1}\sum^n_{r=1}\left(\,\shave{\sum^n_{q=1}}
\frac{\nu\,p_q}{\Omega}\,\frac{\partial}{\partial p_q}\!
\left(\frac{1}{\Omega}\,\frac{\partial H}{\partial x^r}\right)
\frac{\partial H}{\partial x^s}+\frac{\nu}{\Omega^2}\,
\frac{\partial\Omega}{\partial x^s}\,\frac{\partial H}
{\partial x^r}\right)\tau^s_i\,\tau^r_j.
\endgather
$$
Now we are able to write compatibility condition for Pfaff equations
\thetag{8.10}. It breaks into three separate parts. These are the
following equalities:
$$
\gather
\hskip -2em
\gathered
\frac{2\,\nu^2}{\Omega^2}\,\frac{\partial H}{\partial x^s}\,
\frac{\partial H}{\partial p_r}-\nu^2\,\frac{\partial}{\partial x^s}\!
\left(\frac{1}{\Omega}\,\frac{\partial H}{\partial p_r}\right)
+\sum^n_{q=1}\frac{\nu^2\,p_q}{\Omega}\,\frac{\partial}
{\partial p_q}\!\left(\frac{1}{\Omega}\,\frac{\partial H}
{\partial p_r}\right)\frac{\partial H}{\partial x^s}=\\
\vspace{1ex}
=\frac{\nu^2}{\Omega^2}\,\frac{\partial H}{\partial p_r}\,
\frac{\partial H}{\partial x^s}-\nu^2\,\frac{\partial}{\partial p_r}\!
\left(\frac{1}{\Omega}\,\frac{\partial H}{\partial x^s}\right)
+\sum^n_{q=1}\frac{\nu^2\,p_q}{\Omega}\,\frac{\partial}{\partial p_q}\!
\left(\frac{1}{\Omega}\,\frac{\partial H}{\partial x^s}\right)
\frac{\partial H}{\partial p_r},
\endgathered
\tag8.11\\
\vspace{2ex}
\hskip -2em
\aligned
\sum^n_{q=1}\frac{\nu^3\,p_q}{\Omega}\,&\frac{\partial}
{\partial p_q}\!\left(\frac{1}{\Omega}\,\frac{\partial H}
{\partial p_r}\right)\frac{\partial H}{\partial p_s}
+\frac{\nu^3}{\Omega^2}\,\frac{\partial\Omega}{\partial p_s}\,
\frac{\partial H}{\partial p_r}=\\
\vspace{1ex}
&\qquad\qquad=\sum^n_{q=1}\frac{\nu^3\,p_q}{\Omega}\,\frac{\partial}
{\partial p_q}\!\left(\frac{1}{\Omega}\,\frac{\partial H}
{\partial p_s}\right)\frac{\partial H}{\partial p_r}
+\frac{\nu^3}{\Omega^2}\,\frac{\partial\Omega}{\partial p_r}\,
\frac{\partial H}{\partial p_s},
\endaligned
\tag8.12\\
\vspace{2ex}
\hskip -2em
\aligned
\sum^n_{q=1}\frac{\nu\,p_q}{\Omega}\,&\frac{\partial}
{\partial p_q}\!\left(\frac{1}{\Omega}\,\frac{\partial H}
{\partial x^r}\right)\frac{\partial H}{\partial x^s}+\frac{\nu}
{\Omega^2}\,\frac{\partial\Omega}{\partial x^s}\,\frac{\partial H}
{\partial x^r}=\\
&\qquad\qquad=\sum^n_{q=1}\frac{\nu\,p_q}{\Omega}\,\frac{\partial}
{\partial p_q}\!\left(\frac{1}{\Omega}\,\frac{\partial H}
{\partial x^s}\right)\frac{\partial H}{\partial x^r}+\frac{\nu}
{\Omega^2}\,\frac{\partial\Omega}{\partial x^r}\,\frac{\partial H}
{\partial x^s}.
\endaligned
\tag8.13
\endgather
$$
It is easy to check that all these three equalities \thetag{8.11},
\thetag{8.12}, and \thetag{8.13} turn to identities if we substitute
\thetag{5.3} for $\Omega$. Therefore we can formulate the following
main result of this section.
\proclaim{Theorem 8.1} Each modified Lagrangian dynamical system
\thetag{7.1} satisfies {\bf additional normality condition} with
respect to deviation functions \thetag{7.6} determined by its own
\pagebreak Lagrange function $L$.
\endproclaim
    One should note here that Pfaff equations \thetag{8.10} are
not only compatible, but they are also explicitly integrable in form
of equality $H(p,\nu\cdot\bold n(p))=\const$, which is similar to
the equality \thetag{2.11}.
\head
9. Summary and conclusions.
\endhead
   Theorems~7.1 and 8.1 form main result of present paper. Now we
are to understand this result. Thus, we have arbitrary smooth
manifold $M$ without Riemannian metric, but equipped with Lagrange
function $L=L(p,\bold v)$ defining invertible Legendre transformation
$\lambda$ and satisfying two conditions \thetag{7.4} and \thetag{7.5}.
It's clear that these conditions are rather non-restrictive. Despite
to the absence of Riemannian metric, under the above assumptions
\roster
\rosteritemwd=0pt
\item"1)" one can define concept of {\it normal shift} with respect
to geometric structures determined only by Lagrange function $L$;
\item"2)" one can formulate {\it weak and additional normality
conditions};
\item"3)" one can prove that modified Lagrangian dynamical system
\thetag{7.1} satisfies both normality conditions with respect to
geometric structures determined by its own Lagrange function $L$;
\endroster
The results listed above generalize a part of {\it theory of dynamical
systems admitting normal shift} from Riemannian geometry to the
geometry of Lagrangian dynamics. However, this is not complete
generalization. Indeed, in Riemannian case geometry was determined
by metric tensor $\bold g$, while dynamics was determined by force field
$\bold F$. Theory is based on the interplay of these two structures.
Here both geometry and dynamics are determined by Lagrange function
$L$ yet. In further generalizations one should introduce another dynamical
system in $M$ (either Lagrangian or not Lagrangian), and then one should
measure its capability to implement normal shift of hypersurfaces in
geometry determined by $L$.
\head
10. Acknowledgements.
\endhead
     I am grateful to I.~O.~Rasskazov and B.~I.~Suleymanov who communicated
me reference to the book \cite{2}, which I didn't know before.\par
     This research is supported by grant from Russian Fund for Basic
Research (coordinator Ya\.~T.~Sultanaev), and by grant from Academy of
Sciences of the Republic Bashkortostan (coordinator N.~M.~Asadullin). I am
grateful to these organizations for financial support.\par

\enddocument
\end